\documentclass{amsart}
\usepackage[utf8]{inputenc}
\usepackage{amsmath}
\usepackage{amssymb}
\usepackage{amsthm}
\usepackage[a4paper]{geometry}
\usepackage{todonotes}
\usepackage{enumerate}   
\usepackage{cite}
\usepackage{mathtools} 
\usepackage{extarrows}
\usepackage{comment}
\usepackage{hyperref}
\usepackage{graphicx}
\usepackage{placeins}
\usepackage{algorithmic}
\usepackage{subfigure}
\theoremstyle{plain}
\newtheorem{thm}{Theorem}[section]
\newtheorem*{thm*}{Theorem}

\newtheorem{lem}[thm]{Lemma}
\newtheorem{prop}[thm]{Proposition}

\theoremstyle{definition}
\newtheorem{defn}{Definition}[section]
\newtheorem{exmp}{Example}[section]

\theoremstyle{remark}
\newtheorem{rem}{Remark}

\newcommand\inner[2]{\left\langle #1, #2 \right\rangle}

\newcommand{\tco}{\mathcal{T}}
\newcommand{\bo}{\mathcal{L}(L^2)}

\newcommand{\C}{\mathbb{C}}
\newcommand{\R}{\mathbb{R}}
\newcommand{\Rd}{\mathbb{R}^d}
\newcommand{\Rdd}{\mathbb{R}^{2d}}
\newcommand{\Z}{\mathbb{Z}}

\newcommand{\D}{\mathcal{D}}
\newcommand{\I}{\mathcal{I}}

\newcommand{\F}{\mathcal{F}}
\newcommand{\tr}{\mathrm{tr}}
\newcommand{\sfr}{T_\Omega}
\newcommand{\alc}{\text{ALC}}
\DeclareMathOperator*{\argmax}{arg\,max}

\title{Local Structure and effective Dimensionality of Time Series Data Sets}

\author{Monika D\"orfler}
\author{Franz Luef}
\author{Eirik Skrettingland}
\email{monika.doerfler@univie.ac.at, franz.luef@math.ntnu.no, eirik.skrettingland@ntnu.no}

\begin{document}

\maketitle
\begin{abstract} 
The goal of this paper is to develop novel tools for understanding the local structure of systems of functions, e.g.  time-series data points, such as the total correlation function, the Cohen class of the data set, the data operator and the average lack of concentration. The Cohen class of the data operator gives a time-frequency representation of the data set. Furthermore, we show that the von Neumann entropy of the data operator captures local features of the data set and that it is related to the notion of the effective dimensionality. The accumulated Cohen class of the data operator gives us a low-dimensional representation of the data set and we quantify this in terms of the average lack of concentration and the von Neumann entropy of the data operator and an improvement of the Berezin-Lieb inequality using the projection functional of the data augmentation operator. The framework for our approach is provided by quantum harmonic analysis. 
\end{abstract} 
\tableofcontents

\section{Introduction}For high-dimensional, complex data sets, structured dimensionality reduction methods are essential in order
to enable useful further processing\cite{Roweis2000}. The guiding idea behind some approaches is the hypothesis, that successful (machine) learning is made possible by the fact that data of interest live on low-dimensional manifolds as opposed to the dimensionality of the space in which the data are collected a priori. Inspired by this idea, the approach envisaged in this article hinges on the idea that the dimension of data is not at all canonical, but  can vary according to the chosen   representation of the data. Ideally, the choice of representation avoids  loss  of essential information. While image data are directly accessible and rather easy to interpret, time series data such as audio seem to be harder 
to understand, cluster and classify. Time-frequency (TF) methods are often applied to obtain
image-like representations of time-series data such as speech and music~\cite{hinton2012deep}, and to encode the impact of variance over time, \cite{mo19, babrdohasm19, brdohgegrhaklorsc19, dogrbafl18}.  
Applying convolutional neural network (CNN) architectures to the resulting TF-transformed versions of time-series data points has been surprisingly successful in various 
machine learning (ML) tasks. The underlying processes are, however, not entirely understood. One important hypothesis is the assumption that the informative content of the data actually lies on a manifold of 
significantly lower dimensions than their domain. 
It is not clear, however, how these essential parts of data can be made explicit. 

In the work presented in the current article, 
we take a stance toward the identification of 
time-frequency localized components that determine the  entropy of a data set. 
In particular, certain intrinsic structures, which repeat over 
a data set of interest, may be  encoded 
as TF-local components.   TF-local signal components, which repeat over the data points in a data set, probably play a crucial role in training realizations of a given CNN architecture based
on specific data sets~\cite{DiSc14,DaYu12}.  Intuitively, these components determine 
the
coefficients of the convolutional kernels in the lower levels of the network and thus can be expected to carry the essential structure for a certain problem at hand.

While on the one hand, the random  choice of convolution coefficients can lead to improved 
stability to small perturbations~\cite{xu2021robust}  and also give guaranteed reconstruction properties~\cite{ro09}, random initialization of all coefficients 
still is not always a good strategy and may not lead to a favorable architecture realization~\cite{heetal2015delving}. Careful 
investigation of the interaction between structures present in the data and the impact of
convolutions in lower levels of CNNs is thus necessary to better understand the observed phenomena.

In this work, we make a connection between the ubiquitous local averaging of
 TF-representations of  signals in a given data set, 
data augmentation of original time-domain data and their {\it effective dimensionality }\cite{RoyVet07} (ED) via tools from quantum harmonic analysis developed in \cite{Luef:2018conv,LuSk19,LuSk20}.
The key insight guiding these  efforts  is the  association of the data operator $S_\D$ to a system of functions $\D$.   This allows us to capture the structure of the data  and their interaction.  Data augmentation is then 
formalized as the mixed-state localization operator corresponding to $S_\D$. The data operator is the analogue of the density operator in quantum statistical mechanics. This connection leads us to the investigation of various notions of quantum entropies as well as their implications for the structure of the functions in $\D$ and the properties of their interactions. Our results demonstrate that the von Neumann entropy is well-adapted to quantification of the (augmented) data set's entropy via von Neumann entropy of the data density operator $S_\D$. Von Neumann entropy, in turn, is closely related to   effective dimensionality, which was introduced in order to describe the underlying structural properties of a data set as opposed to its a priori dimension.

Our approach uses  tools such as mixed state localization operators and the accumulated Cohen's class of an operator, which have recently been developed in \cite{LuSk19,LuSk20} in work aimed
to connect time-frequency analysis and quantum harmonic 
analysis \cite{Werner}. Furthermore, new  notions, such as  the \textit{total correlation function},  which  captures local  time-frequency correlation between the different data points, and 
an {\it average loss of concentration}, encoding the interaction between data structure and augmentation, are introduced.  In the context of data analysis and the desire to obtain reasonable information by dimensionality reduction, understanding the smoothing action of augmentation, is formalized.

Our manuscript is divided into two parts:
In Part A  we outline the main ideas and results. We illustrate the connections to and implications for data analysis  with 
the help of various examples. In  Part B we put 
forward the theoretical background concerning quantum harmonic analysis, 
von Neumann entropy, an improvement of the Berezin-Lieb inequality for the data augmentation operator and introduce novel notions connecting data analysis and quantum harmonic analysis. Our sharpening of the Berezin-Lieb inequality for the data augmentation operator involves the projection functional that measures how much the data augmentation operator fails of being a projection. In the final statement, the projection functional is replaced by the average lack of concentration which is more suited for our purposes. The proofs of the main results make use of these techniques.  

\section{Part A: Context, Results and Numerical Examples}
\subsection{Concepts and Notation}
\subsubsection{Data sets and data operators}
A data set $\mathcal{D} = \{ f_i\in V , i = 1, \ldots, N\}$ is a collection of  data points $f_i$ in a certain fixed vector space $V$
with innerproduct, whose properties
are not further specified at this point. For two data points $f_i, f_j$, their tensor product $f_i\otimes f_j$
is an operator acting on $V$, defined by 

\begin{equation}\label{Def:TensorPr}
    ( f_i\otimes f_j) (h) = \langle h, f_j\rangle f_i \quad \text{ for } h\in V. 
\end{equation}
We define, for a given data set 
$\mathcal{D}$
the following data operator: 
\begin{equation}\label{Def:DatOp}
    S_\D = \sum_i f_i\otimes f_i.
\end{equation}
As our normalization, we will always assume that $\sum_{i} \|f_i\|_2^2=1$, equivalently that $\tr(S)=1$.
\begin{rem}
If $V= \R^d$, then $S$ is the empirical covariance operator for the multivariate distribution
of the hypothetic underlying density generating the data set. 
While stationary processes are characterized by their spectral density, more complex 
processes may be described by the singular values of their covariance operator, 
which are closely related to the Karhunen-Loeve transform. 
\end{rem}
In data analysis, Karhunen-Loeve transform  is also known as principal component analysis and is defined as follows.
Since $S$ is a self-adjoint operator, it possesses an orthogonal basis of eigenvectors $h_k$,  with corresponding eigenvalues $\lambda_k$, such that 
\begin{equation}\label{Def:DatOpKL}
    S = \sum_k \lambda_k h_k\otimes h_k
\end{equation}

We note here, that the eigenfunctions of $S$ maximize the following expression: 
 %\begin{equation}\label{eq:compeig}
   %   h_k^\Omega = \argmax_{\|f\|_2=1, f \perp h_1^\Omega, \dots h_{k-1}^\Omega} \int_{\Omega} Q_{S}(f)(z) \ dz.
  %\end{equation}
  
  \begin{equation}\label{eq:PCA}
      h_k = \argmax_{\|\psi \|_2=1, \psi  \perp h_1, \dots h_{k-1}} \langle S \psi , \psi \rangle \, ,
  \end{equation}
  which, for the data operator $S$, takes the form 
    \begin{equation}\label{eq:PCA_S}
      h_k = \argmax_{\|\psi \|_2=1,  \psi \perp h_1, \dots h_{k-1}} \sum_i |\langle \psi ,  f_i\rangle |^2
       \end{equation}
 
\subsubsection{Effective dimensionality}
% why rank is not a good measure of complexity - mainly from Vetterli...  
We define the entropy of a given data set $\mathcal{D}$ as the  von Neumann entropy $H_{vN}(S)$ of the associated data operator $S$: 
\begin{equation}\label{eq:VNent}
    H_{vN}(S)=-\tr(S\log (S)).
\end{equation}
The case of maximal correlation (all data points are scalar multiples of a single data point) corresponds to $H_{vN}(S)=0$, thus 
 increasing $H_{vN}(S)$ can be  interpreted as a decrease in the correlations within the data. 

\begin{rem}
\begin{itemize}
    \item In~\cite{RoyVet07},  Roy and Vetterli proposed a concept they denoted as  effective dimensionality, also cf.\cite{Giu20} which is based on the Shannon
entropy of singular values of a given matrix, which may be thought of as 
the collection of data points. Since their concept is closely related to our understanding of entropy of a data set as given in \eqref{eq:VNent}, we adopt their notion of effective dimensionality of $ S = \sum_i f_i\otimes f_i $ as a synonym for the underlying data set's entropy.
    \item The tools of \textit{quantum information theory} allow us to give a more operational interpretation of the statement that an increase in $H_{vN}(S)$ corresponds to a decrease in correlation. Each of the operators $f_i\otimes f_i$ for $i=1,\dots, N$ describes a quantum state, and we consider the following experiment: Alice picks a quantum state $\rho=f_i \otimes f_i$ at random from these states (with equal probability $1/N$ for each state), and sends it to Bob. Bob knows the data set $\mathcal{D}$, but not which particular state $\rho=f_i \otimes f_i$ was sent by Alice. He therefore measures some observable in the state $\rho$ in order to deduce which of the states $f_i\otimes f_i$, $i=1,\dots , N$ he received. Holevo~\cite{Hol73} has shown that an upper bound for the information Bob can deduce about which state he received is given by $H_{vN}(S)$. In other words, $H_{vN}(S)$ quantifies how difficult it is to tell the operators $f_i\otimes f_i$ apart. We consider a data set $\mathcal{D}$ to be correlated if the operators $f_i\otimes f_i$ are difficult to tell apart, in other words if the maximum information $H_{vN}(S)$ accessible to Bob is small. Of course, notions such as accessible information and measurement of an observable can be made precise, and we refer to Chapter 11 of ~\cite{Wil13}.
\end{itemize}
\end{rem}

\subsubsection{Time-frequency analysis, correlation and  local averaging}
For time series data of various origins, it turns out that characteristics which may in principle be captured by Fourier analysis, that is, which 
are appropriately represented by analysing their spectral content, change over time. Furthermore, correlations are observed over both time and frequency; which means that 
components which are localised in time and frequency 
repeat, up to small disturbances and modifications, over the 
entire signal. 
In order to capture these intrinsic structures, 
more often than not, time-series are processed by means of some
methods from time-frequency (TF) analysis before being exploited by machine learning (ML) methods. 
TF methods such as the short-time Fourier (STFT) or wavelet transform yield a two-dimensional representation of a hitherto one-dimensional data point. 

The respective pre-processing is supposed to extract time-frequency structure 
of the signal 
which is believed to be essential to human perception, and introduces invariance to
phase-shifts. Most commonly, the first processing step consists of taking a 
short-time Fourier transform 
\begin{equation}\label{Eq:STFT} 
    V_g f (t, \xi)  = \int_x f(x) \overline{g}(x-t) e^{-2\pi i x\xi}dx 
\end{equation}
followed by a non-linearity in the form of either
absolute value or absolute value squared. 

\begin{rem}
We will subsequently use the notation $\pi (z ) g (x) = g(x-t) e^{2\pi i x\xi}$ for 
$z = (t,\xi)$. The operator $\pi (z )$ is called a time-frequency shift by $z$.
\end{rem}

For a data point $f$ one thus obtains
a first feature stage $F^0$ as follows
\begin{equation}\label{Eg:F0}
    F^0 (z) = | V_g f ( z) |^2 = |\langle f, \pi (z ) g\rangle |^2,
\end{equation}
for some window function $g$.

Due to the structure of the convolutional layers in a CNN, $F^0$ the first, non-linearly generated 
feature, undergoes local weighted averaging in the next and a few subsequent layers. 
Assuming several  weights, also 
called  convolutional kernels, which are all supported in 
a compact set $\Omega$, we obtain the output of the first convolutional layer: 
 \begin{align*}
    F^1 (z,k) &= (F^0\ast m_k )(z) 
    \end{align*}
The consecutive local averaging steps applied in CNNs are crucial in making CNNs  effective for many ML problems. It is the main 
goal of this work to explain the impact of local TF-averaging
on entropy, which 
encodes  the correlation between data points.

\begin{rem}%\todo{Remove or make more concise..
%Why do we believe in time-frequency-local correlation in data}
When considering a
large data set, only functions $h_k$, which are correlated in the sense of 
\eqref{eq:PCA_S} with several or many data points $f_i$ will get the chance to have a significant eigenvalue.
Considering, additionally, the orthogonality condition,
it seems reasonable to assume that the correlation takes place 
in local components,
which repeat over $\mathcal{D}$.  While different expansions are possible, we 
consider an expansion of the data points in a (tight) Gabor frame, 
that is: 
\begin{equation}\label{eq:TF-reprfi}
    f_i = \sum_\lambda c^i_\lambda \pi (\lambda ) g .
\end{equation} The expansion coefficients $c^i_\lambda$ of $f_i$ will be
interpreted as {\it TF-coefficients}, since, for an appropriately localized 
window $g$, they encode the TF-local energy of $f_i$ at a point
$\lambda$ in phase-space.\\
Writing an arbitrary $f$ similarly as 
$f = \sum_\lambda d_\lambda \pi (\lambda ) g$, \eqref{eq:PCA_S}
takes the form
 \begin{equation}\label{eq:PCA_Gabor}
      h_k = \argmax_{\|f\|_2=1, f \perp h_1, \dots h_{k-1}} \sum_i |\sum_\lambda\sum_\mu d_\lambda \langle \pi (\lambda ) g, \pi (\mu ) g\rangle \overline{c^i}_\mu |^2
  \end{equation}
 which may be interpreted, by ignoring the smoothing factor $V_g g$, as finding the 
 vector of TF-coefficients which are maximally correlated, on average and  in a TF-sense,  with the TF-coefficients of the data. 
 \end{rem}
 
 In order to capture the over-all time-frequency correlation structure in a 
data set $\D$, we introduce the following \textit{total correlation function} 
\begin{equation}\label{eq:totcor}
    \widetilde{S}(z):= \sum_{i,j\in\I} |V_{f_i} f_j (z)|^2
\end{equation}
\begin{exmp}\label{Ex1}
%File: TCNnN.m
As an example, we compute the total correlation of a data set of  $500$ signals
with random TF-coefficients and common TF-weight. More precisely, the data points are given as in \eqref{eq:TF-reprfi}, with uniformly distributed coefficients $c^i_\lambda$, to which a common weight 
$w(z) = \sin(2 \pi z)\cdot (1+z)^{-2}$ is applied, so that the  actual data points are given by 
$ f_i = \sum_\lambda c^i_\lambda\cdot w (\lambda )  \pi (\lambda ) g$.\\
Four example data points are shown in
Figure~\ref{Fig:EPTF}, (a), and the data set's total correlation function is plotted in Figure~\ref{Fig:EPTF}, (b). The total correlation captures the TF-behaviour characterized by the weight $w$, as expected.  
\begin{figure}
\centering
 \subfigure[Example Signals]{   \includegraphics[width=0.55\textwidth]{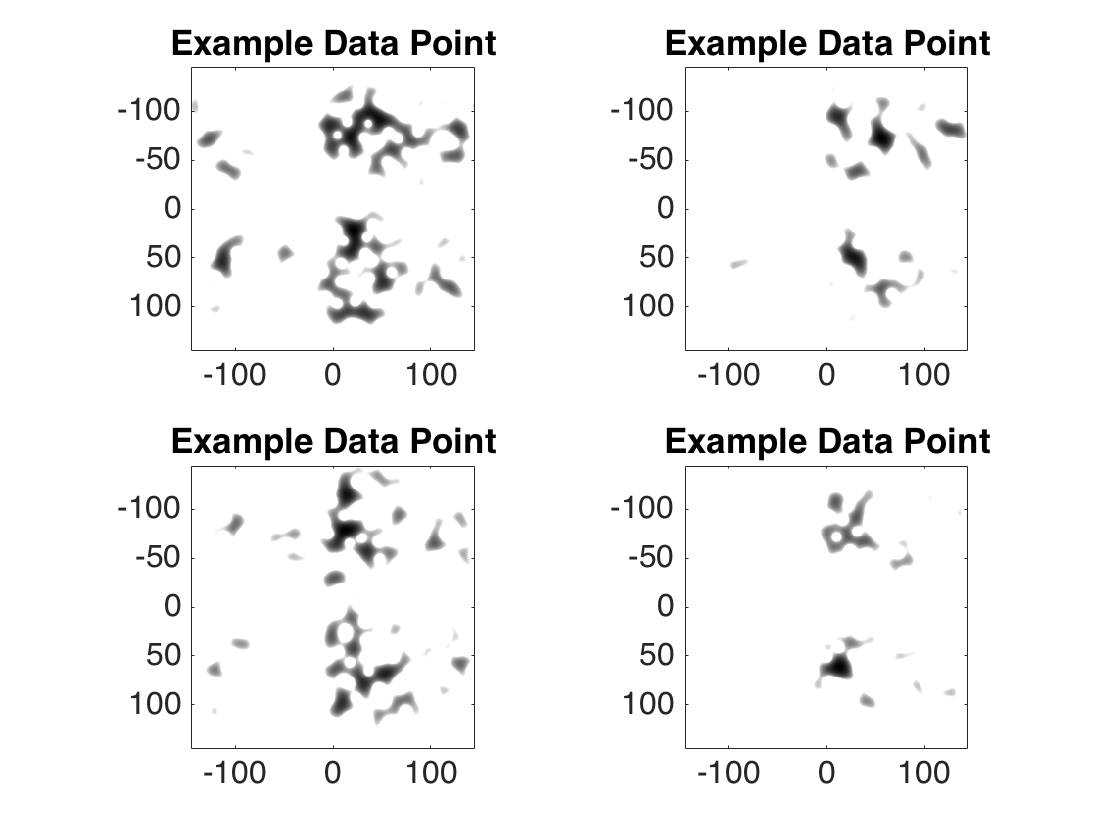}}
   \subfigure[Total Correlation]{     \includegraphics[width=0.35\textwidth]{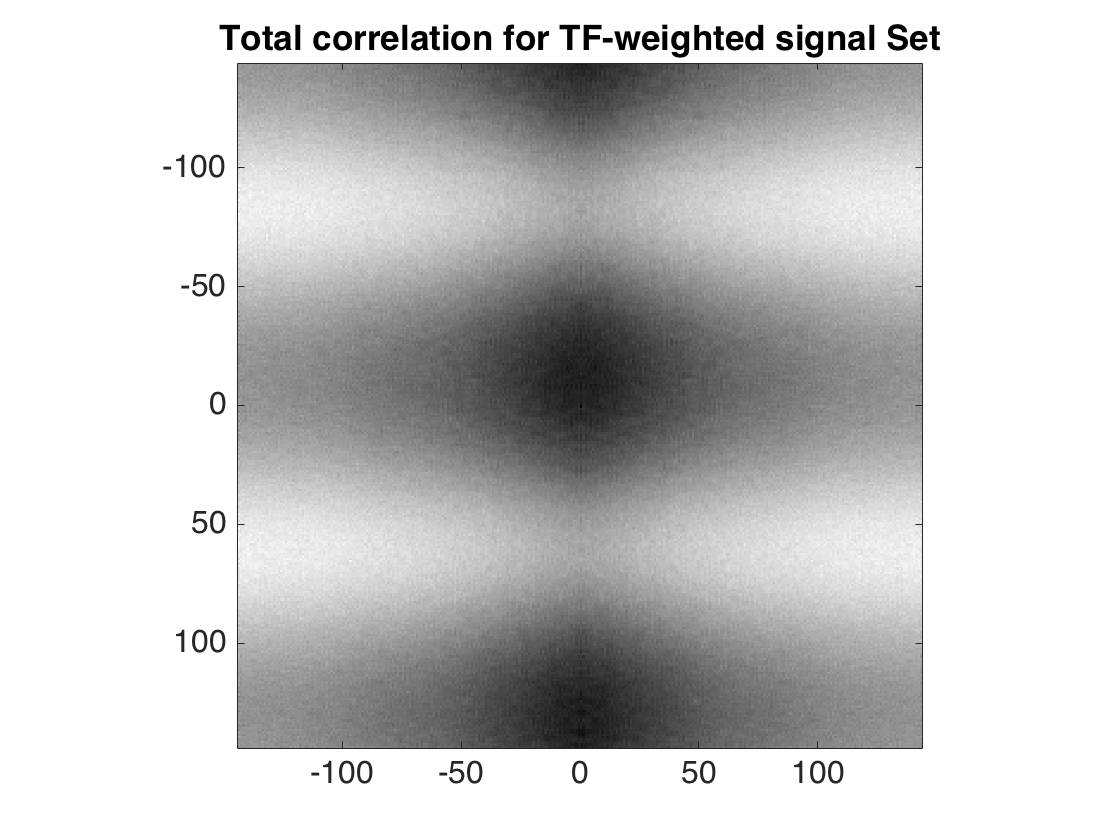}}
\caption{\small\it  (a) Data of random signals with common TF-weight (b) Total correlation $\widetilde{S}(z)$ of data of random signals with common TF-weight.}
    \label{Fig:EPTF}
    \end{figure}

    \end{exmp}

 \subsubsection{Data augmentation} Data augmentation  is a method commonly used
 in data analysis. Its aim is first  to increase the amount of available data by adding carefully modified versions of  data points to existing data, which allows 
 for the usage of larger network architectures. Augmentation 
 often is conducted by  applying certain 
 operators to the elements in $\D$. Second, the particular choice of the operators leads to the enhancement of
desirable invariances of the 
resulting network, whose parameters are learnt from the augmented data set. An example in image processing would be the addition 
of rotated versions of images in the situation of object classification. In the context of audio data, it is reasonable to assume that 
 small TF-shifts have no effect on substantial properties of the data.
 We next give a definition of TF-data augmentation for a  set of TF-sampling points $\mu\in\Omega$. 
 \begin{defn}[$\Omega$-augmentation]
 For a 
 data set $\D$ and a compact set $\Omega\subset\R^2$, its $\Omega$-augmentation is defined by 
 $$\D_\Omega =  \{(|\Omega|^{-1/2}\pi (\mu ) f_i;  f_i\in \D ,\mu\in\Omega\}.$$
 \end{defn}
 The operator corresponding to the augmented data set is then given by
 \begin{equation}\label{eq:TF_augmentation} 
S_{\D_\Omega} := \frac{1}{|\Omega|}\int_{\mu\in\Omega} \sum_i \pi (\mu ) f_i\otimes \pi (\mu ) f_i\ d\mu . %\frac{1}{|\Omega|}\chi_\Omega \star S =
\end{equation}
In Section~\ref{sec:Conv},  we  introduce the concept
of {\it mixed-state localization operators} \cite{LuSk19}. Their definition  
is based on operator convolutions, which will be  formally 
introduced in Definition~\ref{Def:OpConv}. It turns out, that 
 the mixed-state localization operator corresponding  to the data operator $S$
 of a data set $\D$
is precisely the operator corresponding to the $\Omega$-augmented 
 data set $\D_\Omega$,   $\Omega\subset\R^2$. In other words, using the notation 
 from Definition~\ref{Def:OpConv}, in accordance with \eqref{Eq:MSOP} we will also write
  \begin{equation}\label{eq:TF_augmentation1} 
\frac{1}{|\Omega|}\chi_\Omega \star S_{\D}: =S_{\D_\Omega} . 
\end{equation}
\begin{rem}
In machine learning, data augmentation is obviously based on a finite and {\it discrete} set of 
time-frequency shifts, or in fact, any other set of operators which are expected to not influence the data points characteristics in such away that, e.g. class membership changes. 
 \end{rem}
We give some examples in order to illustrate how the mixed-state localization operator corresponding to a data set $\D$ captures information on the data set itself as well as on the impact
of averaging over a domain $\Omega$. 
\subsubsection{Simple motivating examples for the impact of augmentation}\label{Sec:SimpEx}
\begin{enumerate} 
\item 
{\bf Gaussian and  Hermite functions}\\
We first recall the case of Gaussian window $g$, and note that 
we obtain, with $S = g\otimes g$, and  
$$\chi_\Omega \star S (\psi ) =\sum_{\mu\in\Omega}\langle \psi,  \pi (\mu ) g\rangle \pi (\mu ) g,$$  
the situation of 
classical time-frequency localization operators, which have frequently been studied in 
the literature~\cite{Daubechies1988, feno01, defeno02, cogr03-1, abdo12, bagr15}. 
In this simple situation, the interpretation of correlation with respect to the augmentation
\eqref{eq:TF_augmentation} becomes quite straight-forward and can be understood as follows. The total correlation function is 
simply the spectrogram  
$\tilde{S} = |V_g g(z)|^2$, while the eigenfunctions of $S_{\D_\Omega}$ maximise the expression
 \begin{equation}\label{eq:EIGGa}
      h_k = \argmax_{\|\psi\|_2=1, \psi \perp h_1, \dots h_{k-1}} \sum_{\lambda\in\Omega}|\ V_g \psi (\lambda ) |^2
  \end{equation}
and thus result in an orthonormal system of functions, which are optimally correlated with the energy of 
TF-shifted Gaussians $\pi (\mu ) g$, for $\mu$ inside $\Omega$. The resulting function system approximates the Hermite functions, cf.~\cite{Daubechies1988}.
\item 
{\bf Interpolating two Hermite functions}\\ 
We next investigate the interaction between two distinct states and the entropy of the resulting operator
as well as its TF-augmentation.  Let $g$ and $h$ be two $L^2$-normalized functions and define $S_t=\tfrac{1}{2}((1-t) g\otimes g+t h\otimes h)$, for
$t\in[0,1]$.
We first choose $g$ and $h$ to be the first and tenth Hermite function.
In Figure~\ref{Fig:EntIntp}(a),  the entropies of the operators
$ S_t $, for $t$ running from 0 to 1, that is, for the interpolation between first and tenth Hermite function are shown.  Similar behaviour is observed for interpolation between other orthonormal pairs, which shows, that the entropy does not substantially depend on the TF-localization of the involved functions.\\
However,  the entropy-behaviour of the TF-augmentations
$\chi_\Omega\star S_t $ 
is  different. We  consider the interpolation  between Gaussian  and first Hermite function. The corresponding 
entropies of %$S_{\D_\Omega}^{(t)}$ 
$\chi_\Omega\star S_t $ 
are  depicted in  Figure~\ref{Fig:EntIntp}(b). 
We compare these results to the augmentation of the interpolation between zeroth and tenth Hermite function,  depicted in {\bf Figure~\ref{Fig:EntIntp}}(c).

\begin{figure}
        \centering
       
        \subfigure[Entropy $ S_t $]{\includegraphics[width=0.5\textwidth]{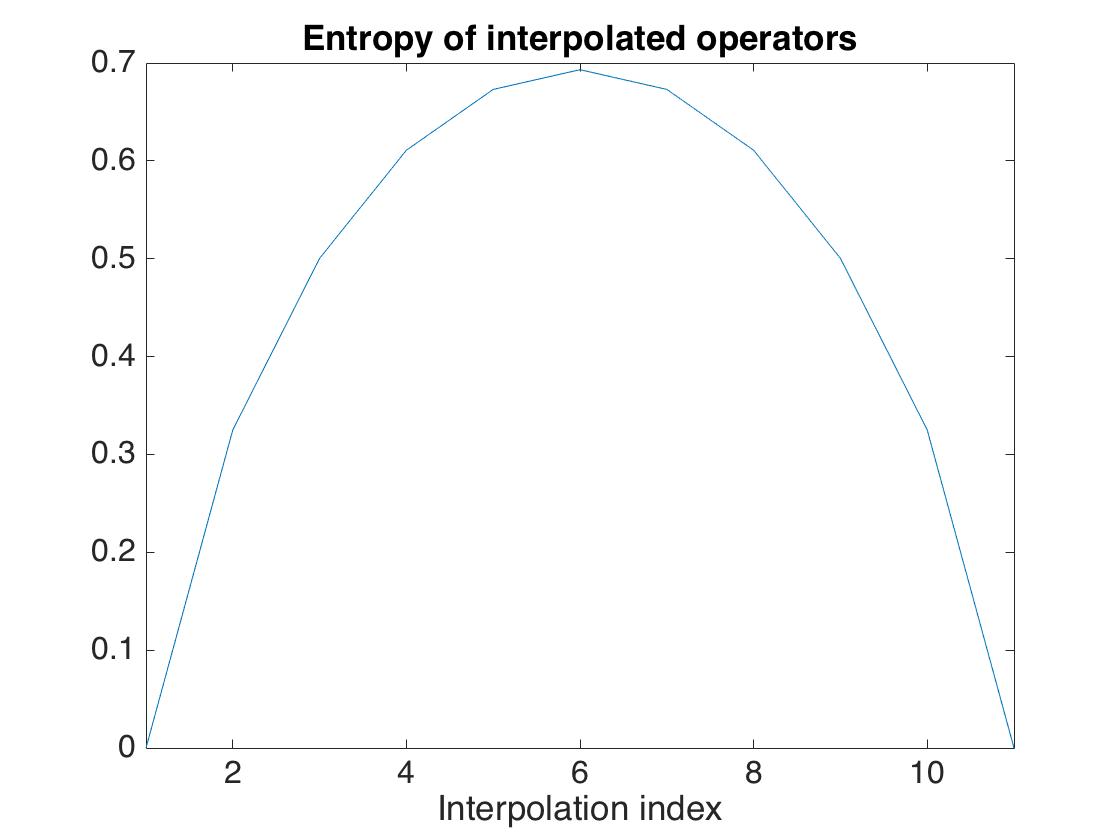}}
        \subfigure[ Entropies of the $\chi_\Omega\star S_t $ for $S_t=\tfrac{1}{2}((1-t) g\otimes g+t h\otimes h)$ with $g$ the first and $h$ second Hermite function.]{\includegraphics[width=0.5\textwidth]{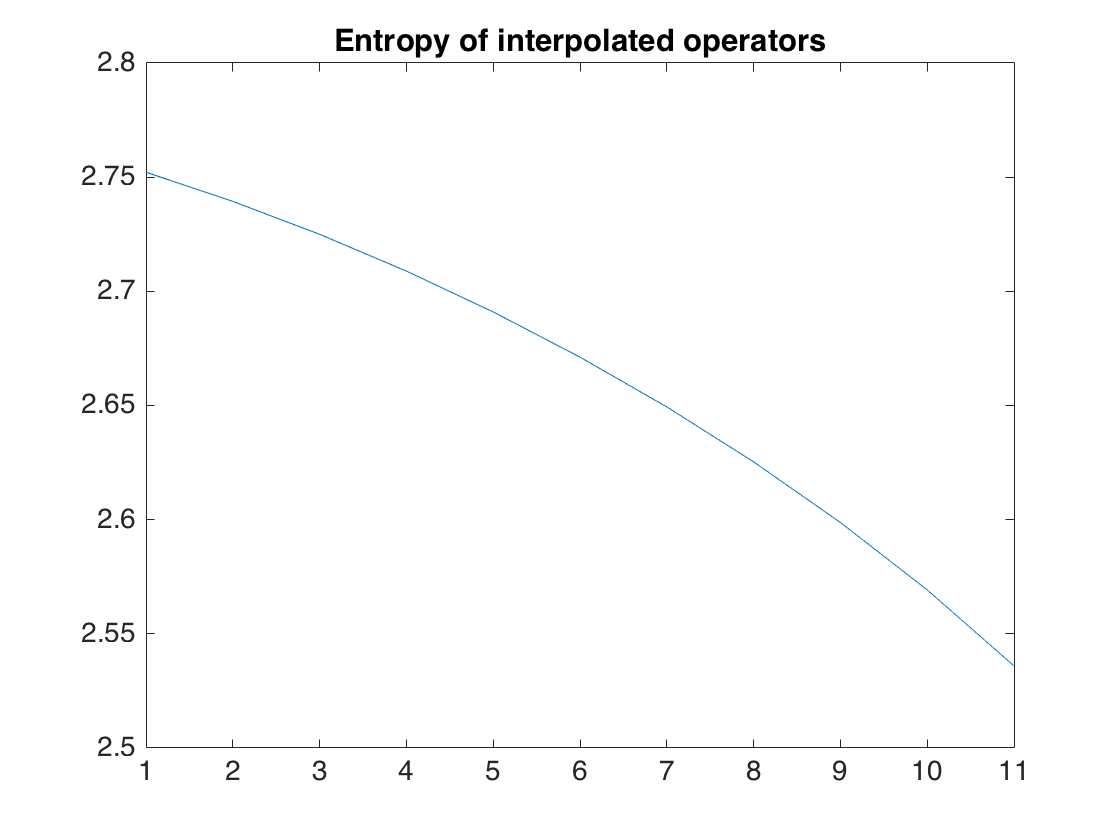}}
         \subfigure[Entropies of the $\chi_\Omega\star S_t $ for $S_t=\tfrac{1}{2}((1-t) g\otimes g+t h\otimes h)$ with $g$ the first and the tenth Hermite function $h$.]{ \includegraphics[width=.5\textwidth]{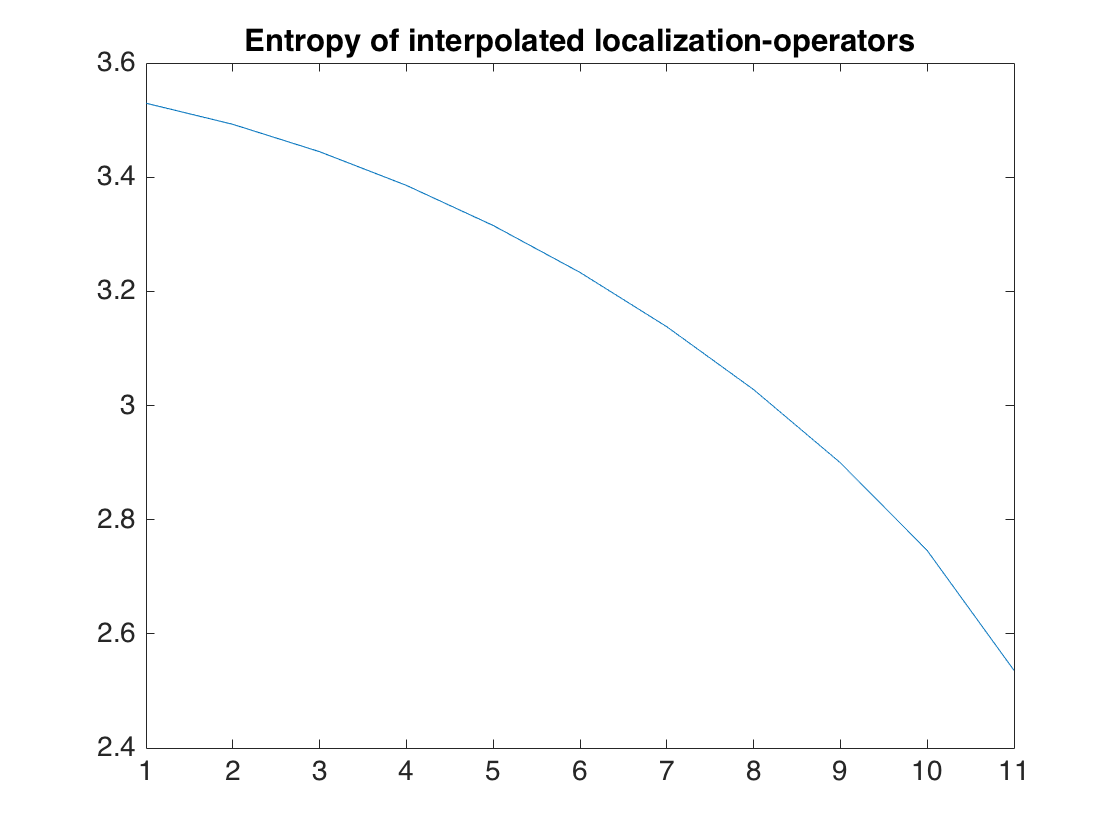}}
         
        \caption{(a) Entropies of the $S_t $, $t$ running from 0 to 1, for interpolation between first and tenth Hermite function. Similar behavior is observed for interpolation between other orthonormal pairs.  (b) Entropies of the $\chi_\Omega\star S_t $ for $S_t=\tfrac{1}{2}((1-t) g\otimes g+t h\otimes h)$ with $g$ the first and $h$ second Hermite function. (c) Entropies of the $\chi_\Omega\star S_t $ for $S_t=\tfrac{1}{2}((1-t) g\otimes g+t h\otimes h)$ with $g$ the first and the tenth Hermite function $h$. }
        \label{Fig:EntIntp}
    \end{figure}
In  Figure~\ref{Fig:SG_Herm}, the TF-behaviour of the three involved functions
is illustrated by their respective spectrogram.  The TF-localization of the single components plays a crucial role for the entropy of the augmented data operator:  the more TF-localized the overall signal energy is within the area of augmentation $\Omega$, the smaller the resulting 
entropy. We next investigate  this insight by means of a more realistic data set.
\begin{figure}  
\includegraphics[width=.8\textwidth]{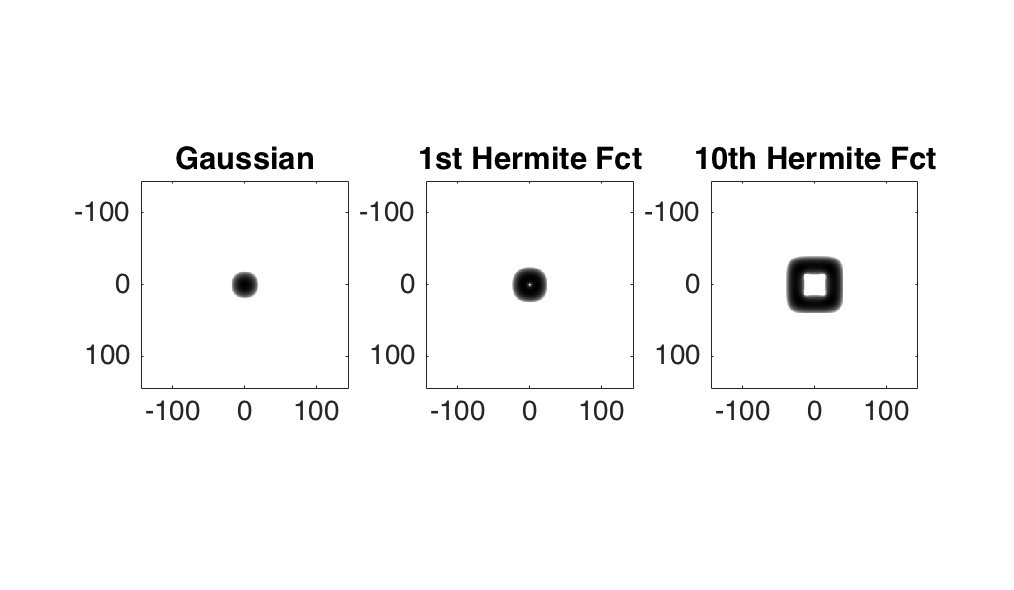}
 \caption{Spectrogram of the zeroth,  first and $h$  tenth Hermite function. }
  \label{Fig:SG_Herm}\end{figure}

\item  {\bf Time-localized chirps}

We construct  a more complex data set with an   inherent TF-
local structure as follows. The 
data set $\mathcal{D}$ consists  of $N$ chirps $f_i$, $i = 1, \ldots , N$, with random base frequency between $30Hz$ and $65Hz$ (normally distributed with mean $50Hz$)  and envelopes given by randomly rotated 
Bartlett-Hanning windows, each of length $280$  in samples; see Figure~\ref{fig:EV_ch}(a) for some examples. 
We let $N$ vary between $100$ and $400$ and compute the rank and effective dimensionality of the corresponding 
data matrix $S_{ch} = \sum_{i = 1}^{280} f_i\otimes f_i$. Rank grows linearly with the number of data points, saturates at the size of the matrix, that is, $280$, and does not reflect the structure of the data set.  
Effective dimensionality, on the other hand, depicted in
Figure~\ref{fig:EV_ch}(b) for several random realizations,
remains relatively stable with respect to the size of the 
data set; the interpretation of this behaviour is that
this measure, i.e. entropy, encodes the actual information
present in the data, which can be extracted given
a sufficient number of data points.\\
The entropies of  the augmented data sets $\D_\Omega$ for 
$\Omega$ a square of side-length 80 are shown by the dotted lines in 
Figure~\ref{fig:EV_ch}(b). The 
effective dimensionality of this data set can be extracted 
from  fewer data points in the original data set $\D$.
On the other hand,  augmentation with respect to this particular 
$\Omega$ significantly increases the data set's entropy. The results 
in this article show that augmentation can be designed as to reduce the increase of entropy. The design of the augmentation depends on the correlation structure in the data. 
To see this, we also plot the total correlation $\widetilde{S_\D}$ for the chirp data set, in Figure~\ref{fig:CHTC}. It is visible, that due to the restricted frequency band, correlation in frequency direction decays quickly, while, given to the uniformly
distributed time-localization, correlation in time is more persistent, since also localized chirps, which
are significantly separated in time can be highly correlated if they have sufficient overlap in frequency. 
\begin{figure}
        \centering
        \subfigure[Data points]{\includegraphics[width=0.39\textwidth]{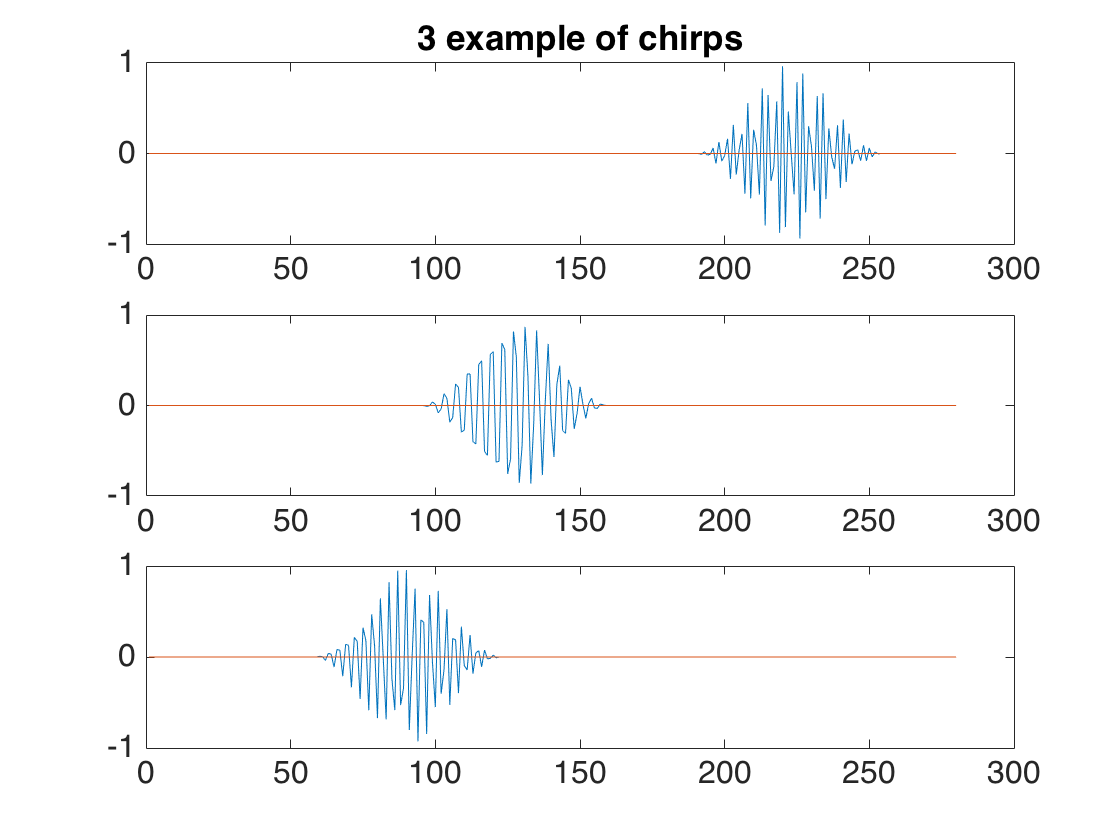}}
        \subfigure[Entropies]{\includegraphics[width=0.59\textwidth]{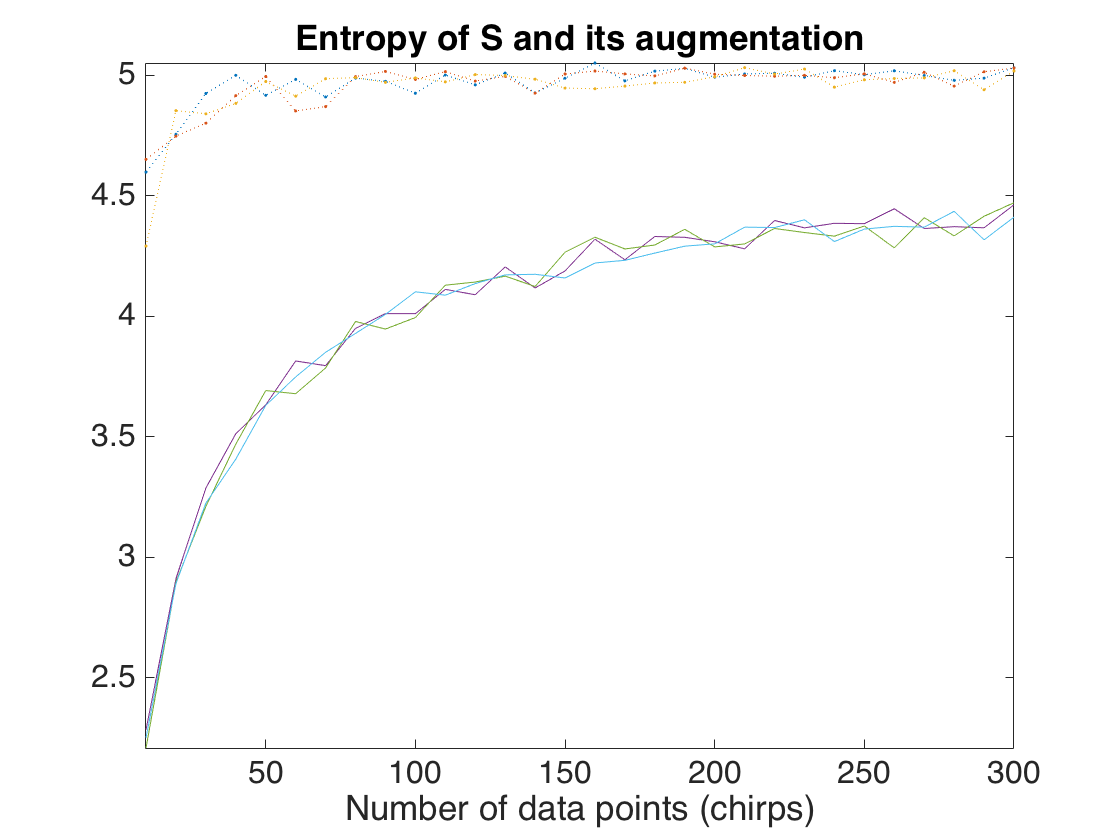}}
        \caption{(a) Windowed chirps (b) Effective dimensionality of data operator for the chirp data set and its $\Omega$-augmentation. Entropy of the augmented data set is dotted. Note that for the latter, the data operator has full rank even for a single data point.  }
        \label{fig:EV_ch}
    \end{figure}

\begin{figure}
        \centering
        \includegraphics[width=0.7\textwidth]{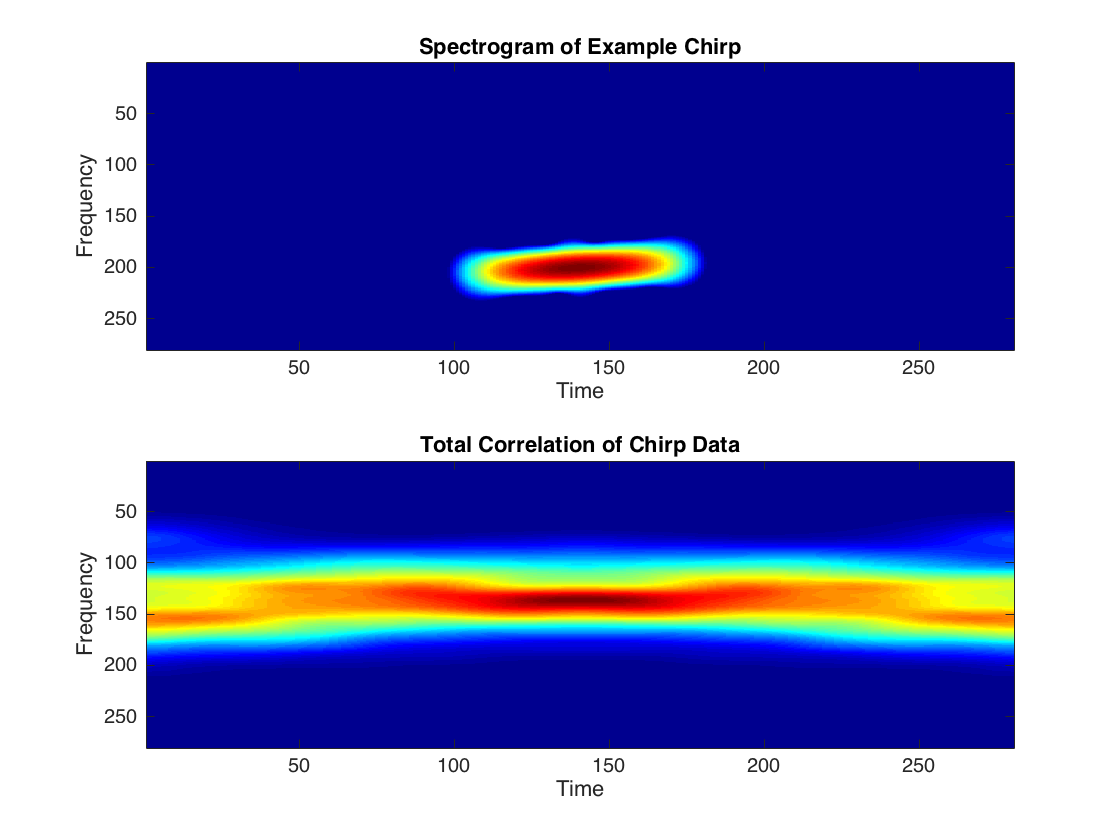}
        \caption{Total Correlation $\tilde{S}$ for Data set of windowed chirps}
        \label{fig:CHTC}
    \end{figure}

\end{enumerate}

In the interpretation of the local correlation in data points, we now deal with the following hypothesis: 
if $\D$ is highly locally correlated in $\Omega$ (e.g. weakly TF-shifted Gaussian windows), then augmentation 
by $\chi_\Omega$ will not increase its entropy significantly, and a small number of eigenfunctions 
of $S_{\D_\Omega}$ represents the augmented  data set well. Our main results, which we present next, give  insight into the relation between a data set, its augmentation by $\chi_\Omega$ and the domain $\Omega$. 
\subsection{Main Results} \label{Sec:MainR}
%\todo{Should we say something about our hopes that low effective dimensionality of augmented data set means local correlations too?}
The effective dimensionality of the augmented data set, represented by the operator $\chi_\Omega \star S$, depends on the original data set represented by $S$, on the domain $\Omega$, and crucially on the interaction between $S$ and $\Omega$. This  interaction reflects the presence of local correlations in the data set, that is, the presence of signal structures which, up to small TF-shifts, repeat over $\D$. \\
Understanding the relation between the effective dimensionality of $\chi_\Omega \star S$ and the pair $(S,\Omega)$ is an important goal of this paper.  In particular, since the effective dimensionality of the augmented data set is defined rather abstractly as the von Neumann entropy of an operator, we aim to understand this abstract object in terms of certain auxiliary concepts whose interpretations are less obscure.  

One such concept is the already introduced total correlation function \begin{equation*}
    \widetilde{S}(z):= \sum_{i,j\in\I} |V_{f_i} f_j (z)|^2
\end{equation*}
when 
\[S=\sum_i f_i\otimes f_i,\]
which encodes the TF-correlations in the data set. In fact,
assume a TF-representation of each $ f_i = \sum_\lambda c^i_\lambda \pi (\lambda ) g $ 
as suggested in 
\eqref{eq:TF-reprfi} and write 
$V_g g ( \mu -\lambda ) = \langle \pi(\lambda )  g, \pi(\mu ) ) g\rangle$, then 
\begin{equation*}
    \widetilde{S}(z):= \sum_{i,j\in\I} | \langle c^i, T_z (c^j\ast V_g g) \rangle|^2.  
\end{equation*} This, by assuming $c^j\ast V_g g\approx c^j$, shows that 
$\widetilde{S}(z)$ encodes {\it TF-local }
correlations for small $z$. The slower $\widetilde{S}(z)$ decays, the 
more TF-distant correlations must be expected over the data set.\\
%\todo{ok formulation, or add something about time-frequency?}. 
Of course, $\tilde{S}$ only gives information regarding the (non-augmented) 
data set
described by $S$. To describe its interactions with $\Omega$,
introduced by the augmentation, we, therefore, propose a new concept called 
    the \textit{average lack of concentration},
defined by 

\begin{equation*}
    \alc(\tilde{S},\Omega) := \frac{1}{|\Omega|} \int_{\Omega} \left( 1- \int_{\Omega-z} \tilde{S}(z') \, dz' \right) \, dz.
\end{equation*}

It turns out that these quantities allow us to give lower and upper bounds for the abstract quantity of effective dimensionality of the augmented data set. In the next theorem, which is Theorem \ref{thm:ALC_ED}, $H(f)$ denotes \textit{differential entropy} of a probability distribution $f$ on $\Rdd$ defined by \[H(f)=-\int_{\Rdd} f(z) \log f(z) \, dz.\] 

\begin{thm*}
For $\Omega\subset \Rdd$ compact and $S$ a positive trace class operator with $\tr(S)=1$,  the following inequalities hold: 
\begin{equation*}
    \log |\Omega|+ \alc(\tilde{S},\Omega) \leq  H_{vN}\left( \frac{\chi_\Omega}{|\Omega|}\star S \right) \leq H\left(\frac{\chi_\Omega}{|\Omega|} \ast \tilde{S} \right).
\end{equation*}
\end{thm*}

\begin{comment}
\subsubsection{Data entropy and local correlations}
In Theorem~\ref{Th1}, we state that the entropy of the data set $\mathcal{D}$ is determined by the extend of local correlation
within the data points in $\mathcal{D}$. More precisely, we use the following definition: 
\begin{defn}[Local data correlation]\label{Def:LocCor}
Fix $\Omega \subset \Rdd$ and a (normalized) data set  $\D = \{ f_i, i = 1,\ldots , N\}$ with associated operator $S=\sum_{i=1}^N f_i\otimes f_i$.
Given $\epsilon >0$, the data set $\D$ is called (internally) $\epsilon-$correlated with respect to $\Omega$ if
\begin{equation}
 \frac{H_{vN}(\frac{1}{|\Omega|}\chi_\Omega \star S)}{\log |\Omega|}<1+\epsilon.
\end{equation}
\end{defn}
We give an informal version of our first central statement. 
\begin{thm}~\label{Th1}
A data set $\D$ is $\epsilon-$correlated with respect to $\Omega$ if its $\Omega$-augmentation has von Neumann 
entropy less than $\log |\Omega| (1+\epsilon)$. This is true, if for the eigenvalues $\lambda_k$ of $ \chi_\Omega \star S$, the following 
inequality holds: 
\begin{equation}
 \sum_k \lambda_k \log (\lambda_k)  > - \epsilon\cdot |\Omega | \cdot \log(|\Omega |).
\end{equation}
\end{thm}
\end{comment}

\subsubsection{Approximation of augmented data operator}
We will also give an approximation of the augmented data operator $\chi_\Omega\star S$  by its eigenvectors,
and approximation quality is quantified by the average lack of concentration. 
As we will see,  for a compact domain $\Omega$ and a density operator $S$ the mixed-state localization operator $\chi_\Omega \star S$ has a diagonalization
\begin{equation*}
    \chi_\Omega \star S =\sum_{k=1}^\infty \lambda_k^\Omega h_k^\Omega \otimes h_k^\Omega, 
\end{equation*}
where $\lambda_k^\Omega$ are the eigenvalues of $\chi_\Omega \star S$ with eigenfunctions $h_k^\Omega$ and $\sum_{k=1}^\infty \lambda_k^\Omega = \tr(\chi_\Omega \star S)=|\Omega|.$

We further let 
$A_\Omega = \lceil |\Omega| \rceil$
and define 
\begin{equation*}
    \sfr = \sum_{n=1}^{A_\Omega} h_k^\Omega \otimes h_k^\Omega , 
\end{equation*}
then we may approximate $\chi_\Omega \star S$ by $\sfr$. The approximation quality depends on
the correlation properties of $\D$, as measured by $\tilde{S}$, and the interaction of correlations and $\Omega$, as measured by the average lack of concentration. In the last inequality we decouple the influence of the correlations and the domain to obtain  an upper bound in terms of the size $|\partial \Omega|$ of the perimeter of $\Omega$.

\begin{thm*}
For $\Omega\subset \Rdd$ compact and $S$ a positive trace class operator with $\tr(S)=1$, let $T_\Omega=\sum_{k=1}^{A_\Omega} h_k^\Omega\otimes h_k^\Omega$ where $A_\Omega=\lceil |\Omega| \rceil$. Then
\begin{align*}
    \frac{\|\chi_\Omega\star S-T_\Omega\|_{\tco}}{|\Omega|} &\leq \frac{A_\Omega-|\Omega|}{|\Omega|}+ 2\cdot \alc(\tilde{S},\Omega) \\
    &\leq \frac{A_\Omega-|\Omega|}{|\Omega|}+ 2\frac{ |\partial \Omega|}{|\Omega|}  \int_{\Rdd} \tilde{S}(z) |z| \ dz.
\end{align*}
\end{thm*}
As our first theorem shows that we can bound $ALC(\tilde{S},\Omega)$ from above using $H_{vN}\left( \frac{\chi_\Omega}{|\Omega|}\star S \right)$, we can easily get an upper bound in terms of $H_{vN}\left( \frac{\chi_\Omega}{|\Omega|}\star S \right)$ as well. We observe that the quality of approximation of the augmented data set by 
its principal components is entirely determined by the interaction between the 
original data set's total correlation and 
the augmentation characteristics. Several examples will enlighten the connection in the
following section. \newpage
\subsection{Examples and Numerical Illustrations}

\subsubsection{Simple motivating examples revisited}
We return to the examples from Section~\ref{Sec:SimpEx} in order to interpret the results on 
effective dimensionality and ALC for the simple data sets presented there. 

\begin{enumerate} 
\item Gaussian functions  and their linear combinations  \\
We first consider  a set of data points generated by a small family of time-frequency shifted Gaussian windows $\pi (\lambda) g$ as follows
\[ f_i = \sum_{l = 1}^{N} c^i_l \pi (\lambda_l^i) g, \, l\in\Z\times\Z\cap M , i\in\I\]
with random coefficients $c^i_l\in\C$ and time-frequency coordinates 
$\lambda^i_l$ randomly chosen in a rectangular domain $M$, 
which is symmetric about $0$ in both time and frequency 
direction, more precisely, of size $2.1875\times 0.3125$. \\
We then compare the average lack of concentration  to the effective dimensionality of the   data set generated by 
augmentation of $\mathcal{D} = \{f_i,\,i\in\I\} $ with respect to  different  domains $\Omega_j$, where 
shape and size of  $\Omega_j$ is allowed to vary during the experiment. More precisely,
$\Omega_1$ is quadratic ($2.45\times 2.45$), $\Omega_2$ is a a wide rectangle of size
 $4\times 1.49$ and $\Omega_3$ is a narrow rectangle of size
 $1.49\times 4$; each of the domains has thus initially size approximately $6$ 
 and  whose side-lengths  enlarged in a second and third step 
 by a factor of $1.3$ and $1.6$, respectively. This yields
 yield  sizes of approximately $10$ and $15$, respectively. The results of $100$ trials are shown
in Figure~\ref{fig:ALC_gauss}, where the mean of each experimental setup is shown together with the variance.  It is obvious that the effective dimensionality grows with the size of the 
augmentation domain. On the other hand, $\Omega_1$, whose shape 
is adapted to the correlation structure of the data set, leads
to lower entropies than the other two choices of $\Omega_j$. 
This effect is also reflected in the ALCs. 
The relation between ALC and ED is further studied in the following example.

\begin{figure}[htbp]
\begin{minipage}{\textwidth}
 \includegraphics[width=\textwidth]{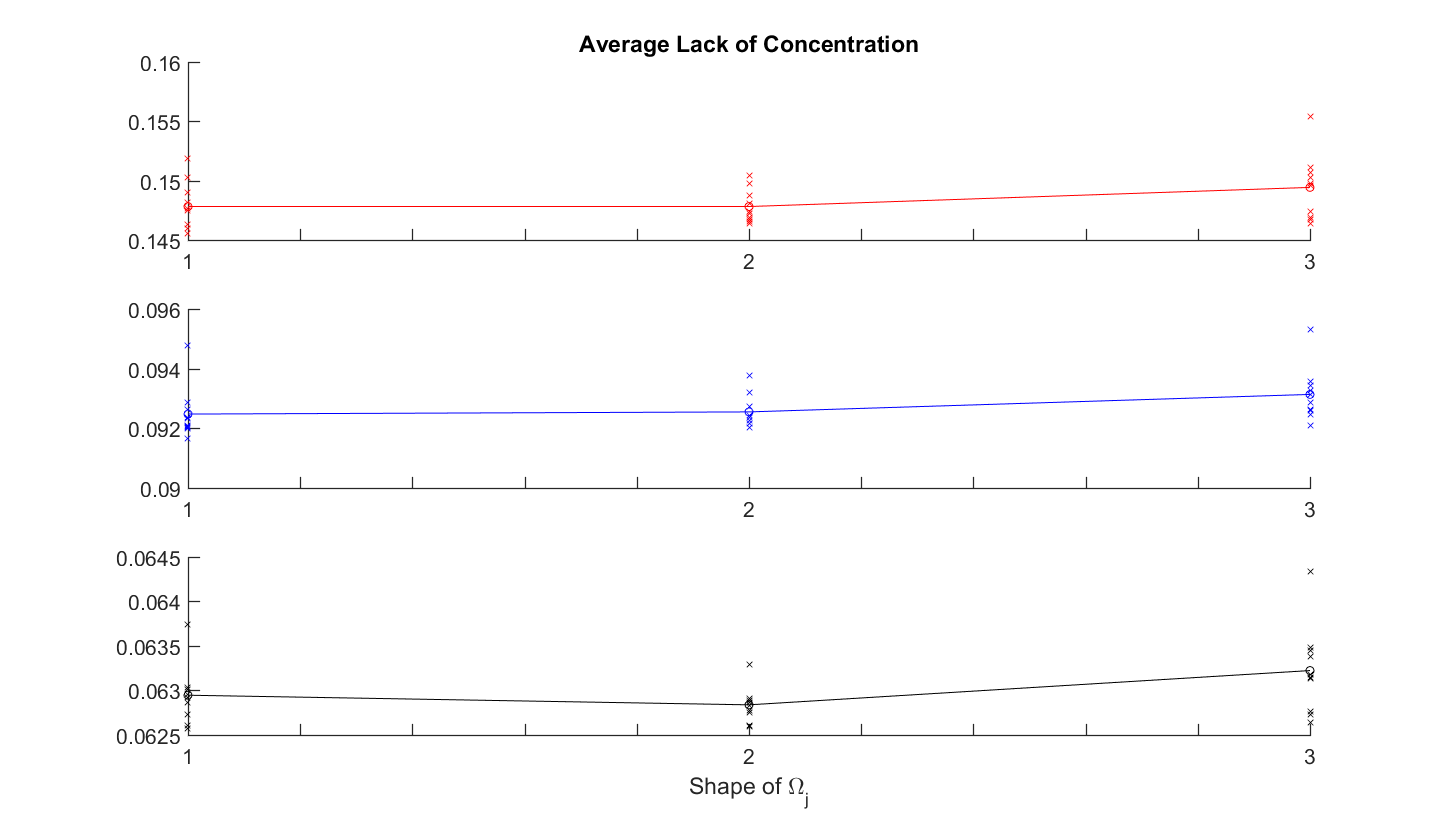}
 \includegraphics[width=\textwidth]{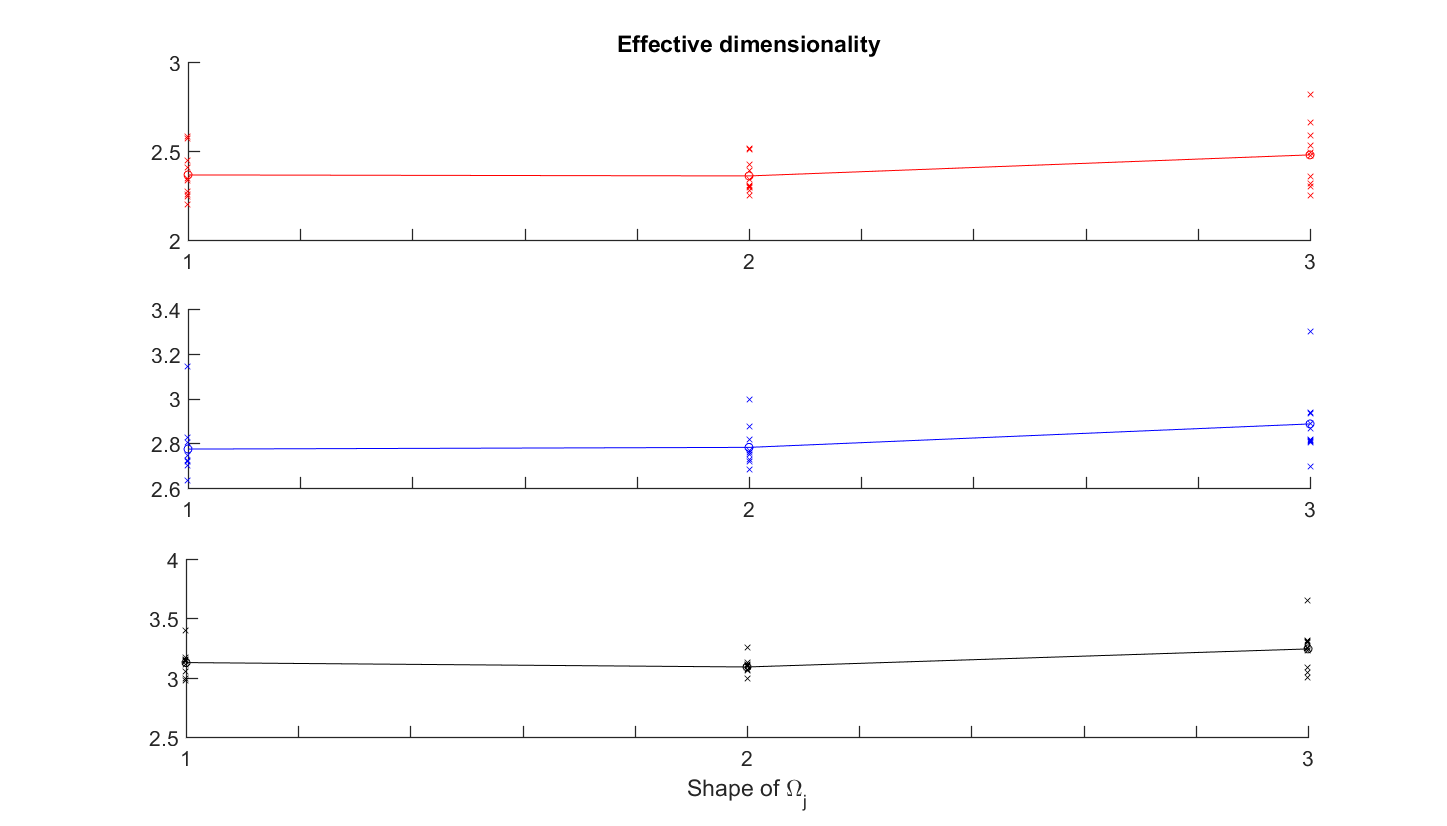}
      \end{minipage}
  \caption{Effective Dimensionality and ALC for data set linear combinations
  of shifted Gaussian windows and varying shape and size of the localization domains $\Omega$, $j = 1,2, 3$: $\Omega_1$ is quadratic ($2.45\times 2.45$), $\Omega_2$ is a narrow rectangle of size
 $1.49\times 4$ and $\Omega_3$ is a wide rectangle of size
 $4\times 1.49$; each of the domains has thus size approximately $6$, upper plots.   Each of  the localization domains  is then enlarged by a factor of $1.3$ and $1.6$,
 respectively, to 
 yield a size of approximately $10$ and $15$, respectively.
 The resulting ALCs and EDs are shown in the middle and lower plots. Each experiment was repeated 100 times, shown are mean (o) and variance (x).}
  \label{fig:ALC_gauss}
       \end{figure}
%%%%%%%%%%%%%CHIRPS
\begin{figure}[htbp]

\begin{minipage}{\textwidth}
 \includegraphics[width=\textwidth]{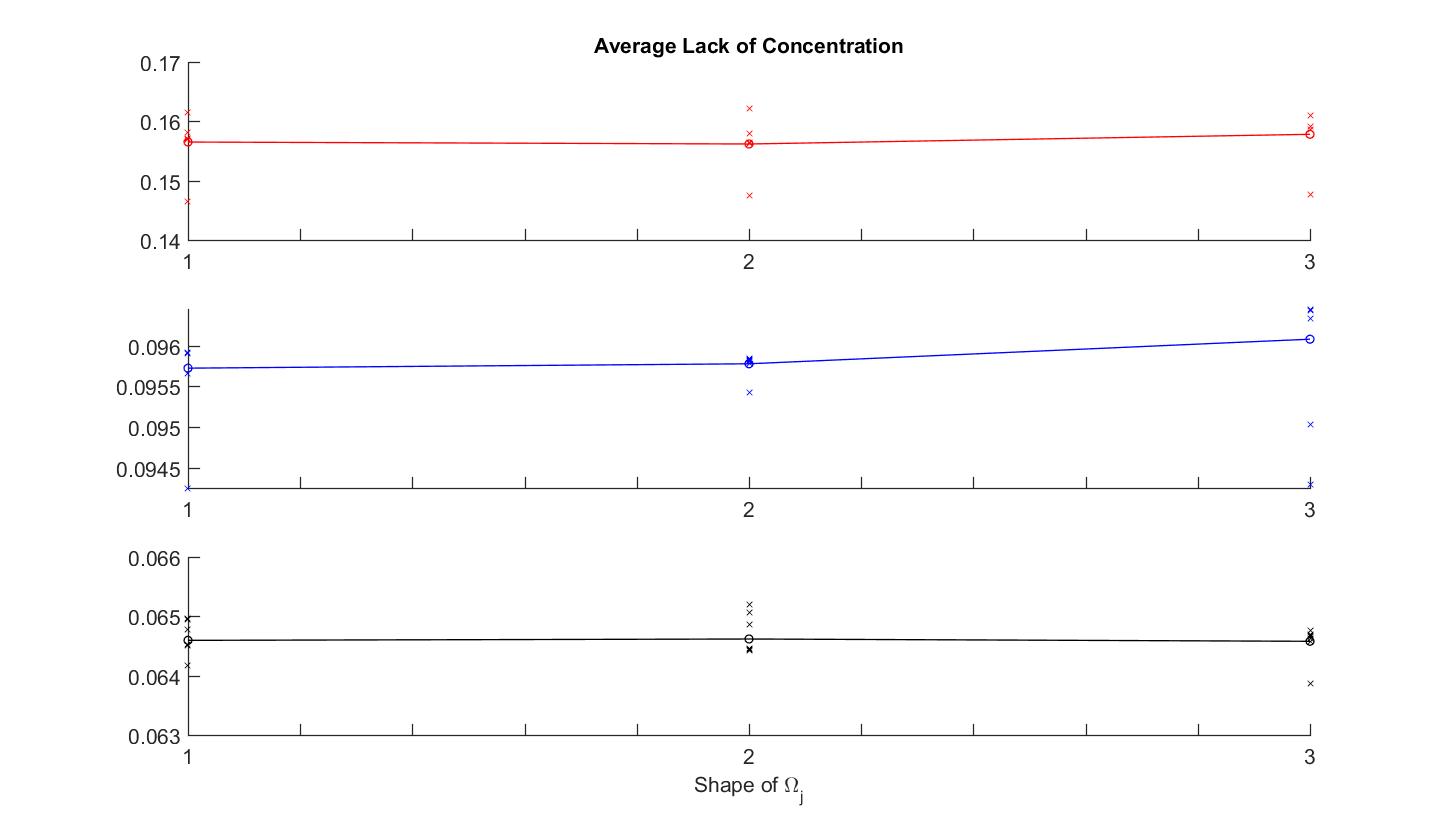}
  \includegraphics[width=\textwidth]{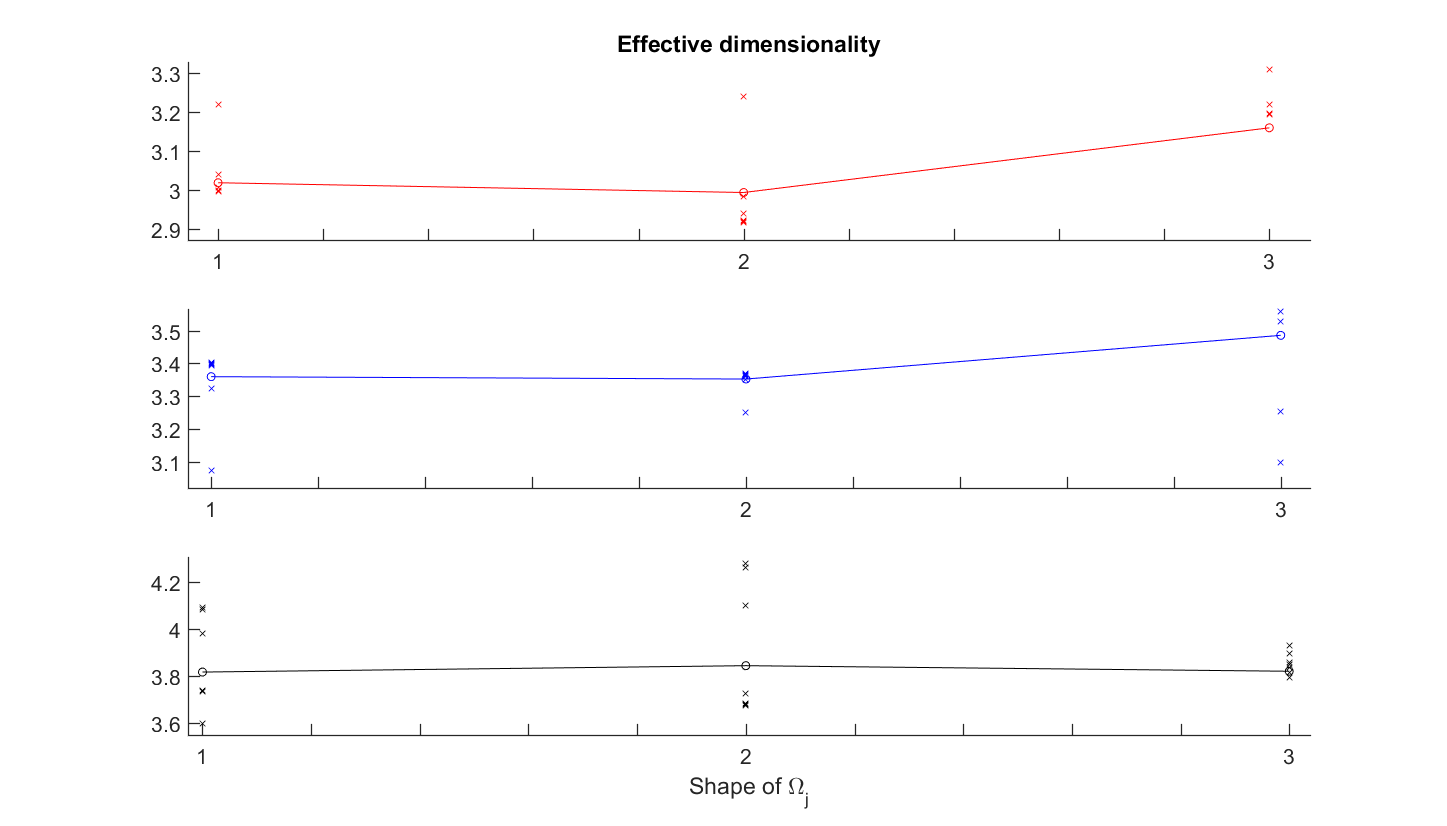}
      \end{minipage}
  \caption{Effective Dimensionality  and ALC for data set of windowed chirps and varying shape and size of the localization domains $\Omega_j$, $j = 1,2, 3$: $\Omega_1$ is quadratic ($2.45\times 2.45$), $\Omega_2$ is a narrow rectangle of size
 $1.49\times 4$ and $\Omega_3$ is a wide rectangle of size
 $4\times 1.49$; each of the domains has thus size approximately $6$.  Each of  the localization domains  is then enlarged by a factor of $1.3$ and $1.6$, respectively, to 
 yield a size of approximately $10$ and $15$, respectively.
  Each experiment was repeated 100 times, shown are mean (o) and variance (x).}
  \label{fig:ALC_Chirpdat}
       \end{figure}

    \begin{figure}[htbp]
    \centering
 \includegraphics[width=1.1\textwidth]{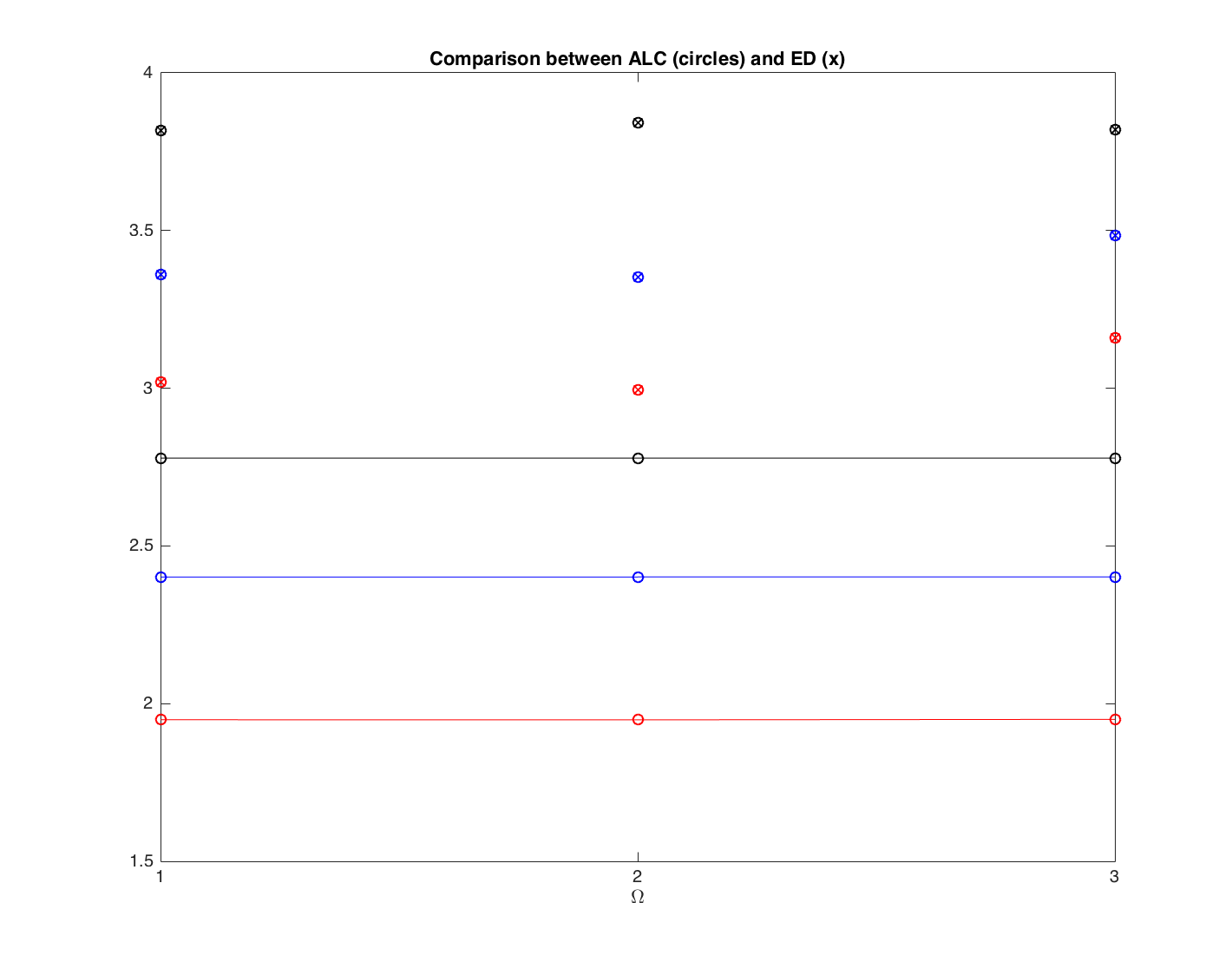}
  \caption{Comparison of Effective Dimensionality 
  ($H_{vN}\left( \frac{\chi_\Omega}{|\Omega|}\star S \right)$ and 
  $  \log |\Omega|+ \alc(\tilde{S},\Omega)$
  for data set of windowed chirps and varying shape and size
  of the localization domains $\Omega$, $j = 1,2, 3$: 
  $\Omega_1$ is quadratic ($2.45\times 2.45$), $\Omega_2$ is a narrow rectangle of size
 $1.49\times 4$ and $\Omega_3$ is a wide rectangle of size
 $4\times 1.49$; each of the domains has thus size approximately $6$ (red values).  
 Each of  the localization domains  is then enlarged by a factor
 of $1.3$ and $1.6$, respectively, to 
 yield a size of approximately $10$ (blue values) and $15$ (black values), respectively.}
  \label{fig:ALC_ED} 
       \end{figure}
 %\item Interpolating Hermite functions \\
%eturning to the observations concerning Hermite functions in Section~\ref{Sec:SimpEx}, we now investigate the 
% total correlation function and the ALC    corresponding to the interpolation between two Hermite functions and their augmentation. 
       
       \item  Time-localized chirps\\
       We study the data set of time-localized chirps, 
       which was introduced in Section~\ref{Sec:SimpEx}. For the subsequent 
       experiments, $150$ data points, generated as 
       described in Section~\ref{Sec:SimpEx}, were used.
 In Figure~\ref{fig:CHTC} we see that total correlation is concentrated in
 a region which is eccentric  with respect to time- and frequency axis. 
 That is, $\tilde{S}$ is better concentrated in frequency 
 direction than in time.  We calculate  
 $\alc(\tilde{S},\Omega)$
 and 
 $H_{vN}\left( \frac{\chi_\Omega}{|\Omega|}\star S \right)$ 
 for both quadratic and  rectangular domains 
 $\Omega$ of varying eccentricity and approximately constant size
 as in Experiment (1) above. \\
 As expected, the shape of $\Omega$ influences ALC and thus 
 effective dimensionality of the augmented data set. Results are 
 shown in Figure~\ref{fig:ALC_Chirpdat} for the three different domains
 $\Omega_j$, $j = 1,2, 3$. 
 Figure~\ref{fig:ALC_ED} illustrates the result
 of Theorem~\ref{thm:ALC_ED}  in more detail, by comparing, 
 for the same experimental setup, the values of 
  $  \log |\Omega|+ \alc(\tilde{S},\Omega)$ and 
  $H_{vN}\left( \frac{\chi_\Omega}{|\Omega|}\star S \right)$.  The illustration of the Theorem's result is obvious; this example makes the relation between the shape of the domain of augmentation, its size $|\Omega|$ and the total correlation of $\D$ more precise. 
 
    \end{enumerate} 
\subsubsection{A data set of local components}
In this example, we consider  a given signal $g\in L^2(\mathbb{R})$, whose over-all energy is assumed well-concentrated within $\Omega$. A  data set 
$\D = \{f_i, \, i \in\I\}$ is generated by 
 data points defined as 
\begin{equation}\label{Eq:ds1}
	f_i=\tilde{f_i}+c_i \pi(z_i)g
\end{equation}
for  $c_i\in\C$ with $0\leq |c_i|^2 \leq 1$ 
and an orthogonal  noise component $\tilde{f_i}\perp \pi(z_i)g$, 
such that $\|f_i\|_2 = 1$. 
 \begin{figure}
        \centering
       \subfigure[Classical Localization Operator]{\includegraphics[width=0.8\textwidth]{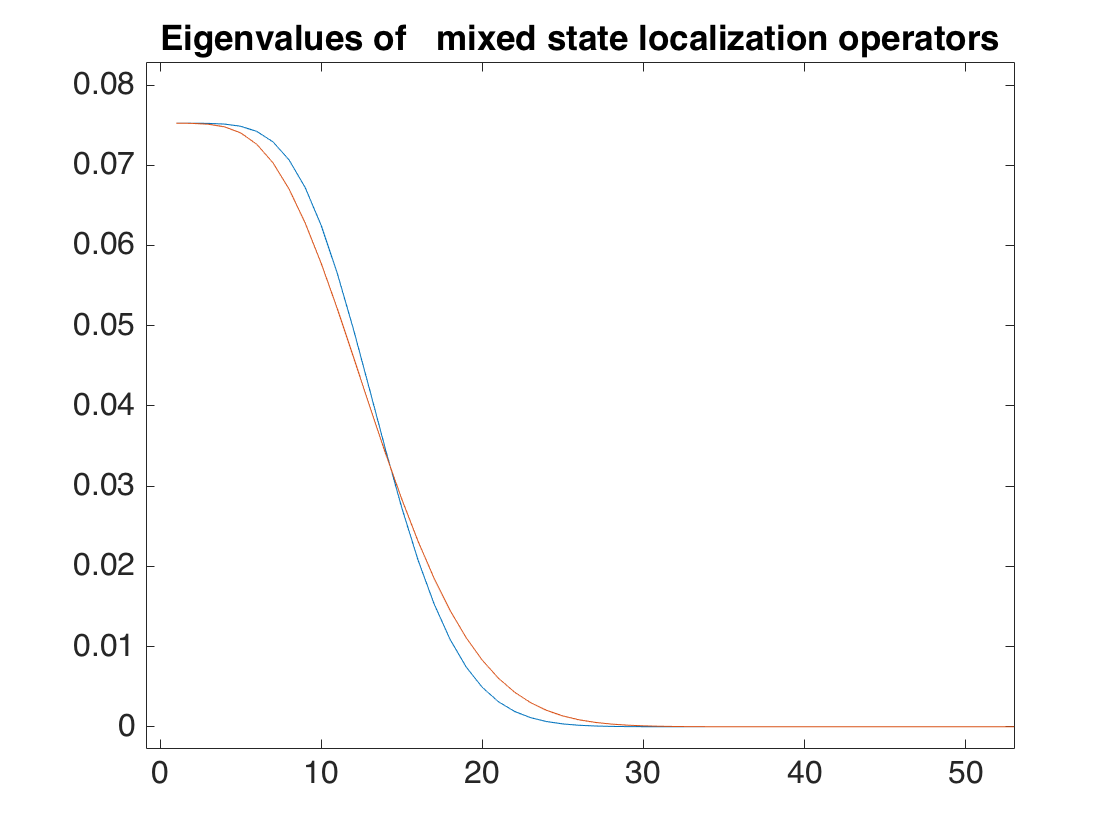}}
        \subfigure[Mixed-state localization operators]{\includegraphics[width=0.8\textwidth]{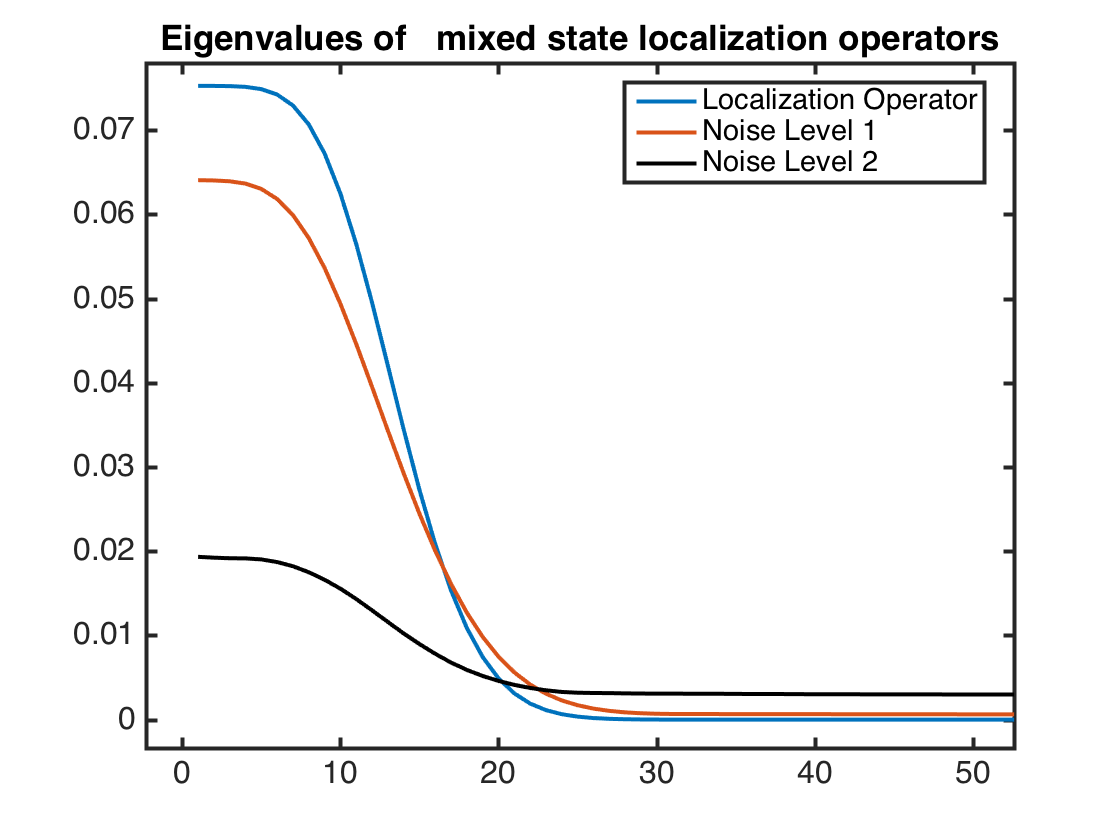}}
        \caption{(a) Eigenvalues of Classical Localization operator (blue) and mixed-state localization operator with 30 Gaussian windows (red). Here, the sequence of coefficients is such that $|c_i| = 1, \,\forall i$ and the uncorrelated noise part is $0$. (b) Eigenvalues of Classical Localization operator (blue) and two versions of  mixed-state localization operators with 30 Gaussian windows. Here, the sequence of coefficients $c_i$ is random and the uncorrelated noise part is  increasing in Level 1 (red) and 2 (black).  }
        \label{Fig:Eigt_Noiseless}
  \end{figure}
 The data set is first generated without noise component and  augmented with respect to a 
    quadratic domain $\Omega$ of size $15$. 
    In Figure~\ref{Fig:Eigt_Noiseless}, the eigenvalues
    of the resulting mixed-state localization operators are shown.
    We see that adding several close Gaussian windows does not significantly change the augmented data
    set's entropy: for the standard localization operator it is $2.86$, while for 
    a data set $\D$ of $30$ TF-shifted Gaussians it is $2.91$, see Figure~\ref{Fig:Eigt_Noiseless}(a). This is intuitively clear, since the overall total correlation remains concentrated inside $\Omega$. The situation starts to change, as soon as the local components explain only part of the data structure, that is, the energy in the noise part $\tilde{f_i}$ of each data point increases, hence the total correlation is not exclusively concentrated inside $\Omega$. The resulting eigenvalues are depicted in see Figure~\ref{Fig:Eigt_Noiseless}(b). Level 1 noise has carries approximately $10\%$ and Level 2 noise $30\%$ of the signal energy.

\subsubsection{The mix of local and global influence}
\begin{figure}
        \centering
        \includegraphics[width=0.7\textwidth]{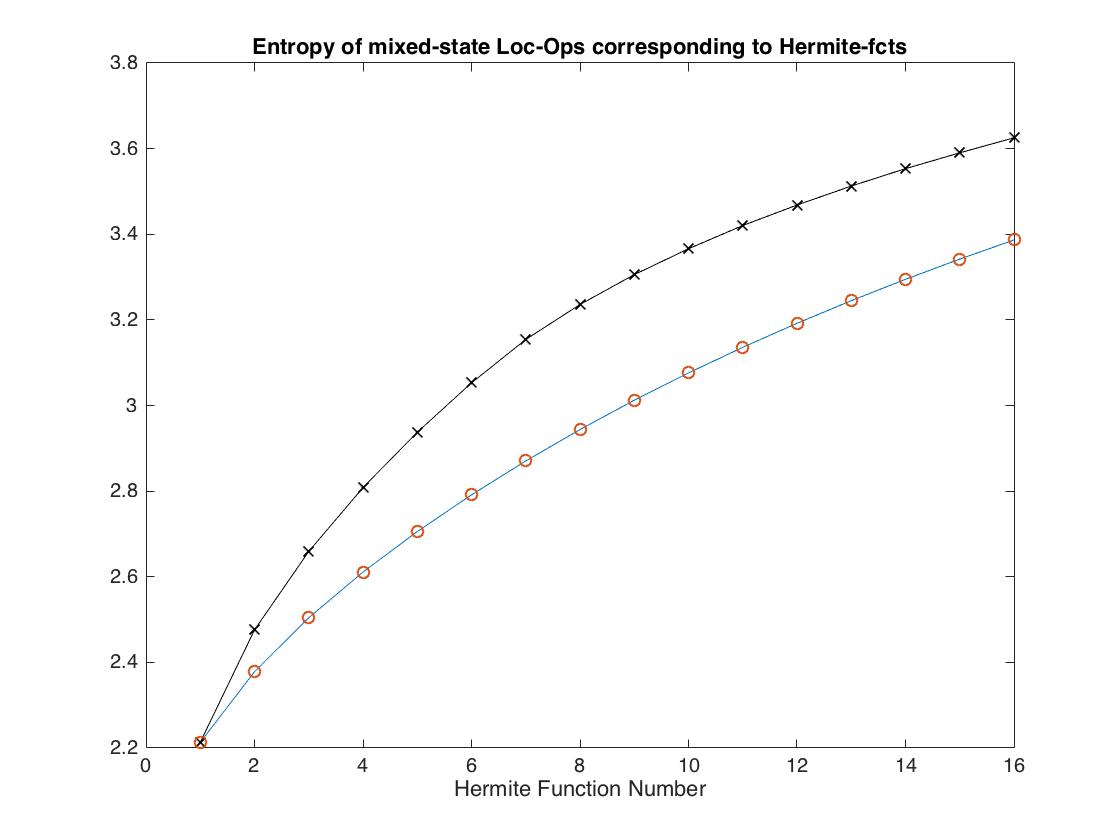}
        \caption{The entropies of mixed-state localization operators corresponding to the projection onto a single 
        Hermite function $h_h$  (black crosses) or onto the set of first $n$ Hermite functions  (red circles).  }
        \label{Fig: EntHermite}
    \end{figure}

This final example illustrates the impact of the TF-characteristics of each data point, 
 the interaction between local
TF-components and their influence on a data set's entropy.  We consider  Hermite functions and distinguish two cases of 
related data sets in order to under 
the influence of the growing TF-spread  of 
Hermite functions $h_i$, $i = 1,...,16$. To this end, we consider the rank-one operators $S_n = h_n \otimes h_n$, which correspond to the data sets $\D_n = \{h_n\}$,  $n = 1,...,16$, and compare them to   $S_1^n = \sum_{i = 1}^n h_i \otimes h_i$, corresponding to 
$\D_1^n = \{h_1,\ldots , h_n\}$. We compute the  entropy
of the augmented data sets $\D_{n,\Omega}$ and $\D_{1,\Omega}^n$, for $\Omega$ a square  of size $9$. 
The second setup, i.e. the augmentation of $\D_1^n$,  actually yields smaller entropy of the
corresponding mixed-state localization operator
$\chi_\Omega\star S_1^n$ than $\chi_\Omega\star S_n$ see Figure~\ref{Fig: EntHermite}. This example shows that augmentation of
smaller data sets may lead to higher effective dimensionality than augmentation of a data set with more data points, even if those are orthogonal. At closer inspection,
this is not a surprising result, since the decay of the eigenvalues of $\chi_\Omega\star S_1^n$ is faster and thus the ALC of $\chi_\Omega\star S_1^n$ is smaller than for $\chi_\Omega\star S_n$. In the context of data augmentation this means that the augmentation is
adapted to the correlation structure of $\D_{1}^n$ rather than to that of $\D_n$. Hence, 
for $\D_n$ augmentation changes the information 
in the data set more significantly than for $\D_{1}^n$, for which the overall, average correlation is better concentrated inside $\Omega$.

\subsection{Interpretation  and perspectives}
In the previous section, we presented several 
data sets of structured time-series data 
in order to illustrate the main results given in Section~\ref{Sec:MainR}. The 
given examples show
an intricate connection between TF-localized correlations in the signals 
contained in a given data set and the effect of augmentation. It is interesting 
to observe, that augmentation, which is used as a means to improve
generalization properties of networks learned from data, leads to varying effects
in dependence on the TF-structure of the original data set. In particular, 
if the structure of the set $\Omega$ is aligned with the area of local 
correlations within the data set, the consequence  of augmentation on the effective
dimensionality is limited. Exploiting this observation, data augmentation can be designed in order to increase the size of a data set without changing its intrinsic structure.  On the other hand, preliminary experiments show the potential of 
adapted augmentation, in order to separate signal classes
with distinct correlation properties. The details and applicability 
of this idea will be elaborated in future work. %Say more?

%%%%%%%%%%%%%%%%%%%%%%%%%%%%%%%%%%%%%%%%%%%%%%%%%%%%%%%%%%%%%%%%%%
\section{Part B: Technical Details and Proofs}
\subsection{Technical Background} \label{sec:techback}
As a theoretical framework for our mathematical results, we will use the theory of quantum harmonic analysis first introduced by Werner in \cite{Werner}. This framework allows us to exploit intuitions and results from harmonic analysis of functions, such as convolutions, to prove results that apply to both operators and functions.
\subsubsection{Convolutions}\label{sec:Conv}
For our purposes, the most important concept from quantum harmonic analysis is that of convolutions of operators and functions. To introduce these operations, we first define a time-frequency shift of an operator as follows. 
\begin{defn}[Shift of operator] Given an operator $S$ on $L^2(\mathbb{R}^d)$, we define its translation $\alpha_z(S)$ to be the operator
\begin{equation*}
  \alpha_z(S)=\pi(z)S\pi(z)^*,
\end{equation*}
where $\pi(z)$ denotes the time-frequency shift of $f\in L^2(\R^d)$ by $z=(x,\omega)\in\Rdd$, $\pi(z)f(t)=e^{2\pi it\omega}f(t-x)$.
\end{defn}
Using this translation operation, we can define two new convolution operations:
\begin{defn}[Operator convolutions]\label{Def:OpConv}
\begin{enumerate}
    \item Let $f\in L^1(\Rdd)$ and let $S$ be a trace class operator on $L^2(\Rd)$. The convolution $f\star S$ is the operator on $L^2(\Rd)$ defined by \begin{equation*}
        f\star S := \int_{\R^{2d}}f(z)\alpha_z(S) \ dz;
    \end{equation*}
    interpreted either as a Bochner integral or in the weak sense:
    \begin{equation*}
	\inner{(f\star S)\psi}{\phi}=\int_{\Rdd} f(z) \inner{\alpha_z(S)\psi}{\phi} \ dz, \hskip 1em \text{ for $\psi,\phi \in L^2(\Rd)$.}
\end{equation*}
\item If $S,T$ are trace class operators on $L^2(\Rd)$, the convolution $S\star T$ 
is the function on $\Rdd$ defined by \begin{equation*}
  S \star T(z) = \tr(S\alpha_z (\check{T})) \hskip 1em \text{for } z\in \Rdd,
\end{equation*}
where $\check{T}=PTP$ for $Pf(x)=f(-x)$.
\end{enumerate}
\end{defn}
We will soon restrict our attention to certain classes of functions and operators, 
where more insight can be gained about these convolutions, but the reader should already note that the convolution of a function with an operator yields a new operator, while the convolution of two operators yields a new function.  

Now assume that $S$ is a positive trace class operator with $\tr(S)=1$ (a so-called density operator), and $f=\chi_\Omega$ for some compact subset $\Omega \subset \Rdd$. The convolution $\frac{\chi_\Omega}{|\Omega|}\star S$ is then the mixed-state localization operator \begin{equation}\label{Eq:MSOP}
    \frac{\chi_\Omega}{|\Omega|}\star S=\frac{1}{|\Omega|} \int_\Omega \alpha_z(S) \, dz
\end{equation}
that we already met when discussing $\Omega-$augmentation in \eqref{eq:TF_augmentation} --- in that case $S$ was a data operator from \eqref{Def:DatOp}. In other words, $\Omega$-augmentation can be described as a convolution operation. 

The total correlation function, which we met in \eqref{eq:totcor} when $S$ was a data operator, also has a simple description as a convolution, namely:
\begin{equation*}
    \tilde{S}(z):=S\star \check{S}(z).
\end{equation*}
Of course, this expression for $\tilde{S}$ makes sense for any trace class operator $S$, but we will mainly consider $\tilde{S}$ when $S$ is a data operator. 

\subsubsection{Properties of the convolutions}
The convolutions introduced above share many properties of the usual convolution of two functions. They are both commutative and associative, and in particular the associativity is a non-trivial and useful property. As an example, we have the relation $$(f\star S)\star T=f\ast (S\star T),$$ where $\ast$ denotes the familiar convolution $f\ast g(z)=\int_{\Rdd} f(z')g(z-z') \ dz'$ of functions. When inspecting this equation we see that it shows the compatibility of three different convolutions: of two functions, of two operators and of a function with an operator.

We summarize some other properties of the convolutions in the following proposition. Proofs can be found in \cite{Werner} or \cite{Luef:2018conv}.
\begin{prop} \label{prop:properties}
Let $S,T\in \tco$.
\begin{enumerate}
    \item (Young's inequality) For any $1\leq p \leq \infty$ and $f\in L^p(\Rdd)$ we have that 
\begin{equation*}
    \|f\star S\|_{\tco} \leq \|f\|_{L^p} \|S\|_{\tco}
\end{equation*}
where $\|\cdot \|_{\tco}$ is the trace class norm. 
\item $S\star T\in L^1(\Rdd)$ with \[\int_{\Rdd} S\star T(z) \ dz=\tr(S)\tr(T).\]
\item $\|S\star T\|_{L^\infty}\leq \|S\|_{\tco} \|T\|_{\bo}$.
\item If $f\in L^1(\Rdd)$, then $\tr(f\star S)=\int_{\Rdd} f(z) \ dz \cdot  \tr(S).$
\item If $S,T$ are positive operators, then $S\star T$ is a positive function.
\item If $f\in L^1(\mathbb{R}^{2d})$ is a positive function and $S$ is a positive operator, then $f\star S$ is a positive operator.  
\end{enumerate}
 
\end{prop}
Versions of the so-called \textit{Berezin-Lieb} inequalities have been shown to hold in quantum harmonic analysis by different authors \cite{Werner,Klauder:2011}. Since none of the existing formulations in the literature cover the case we are interested in explicitly, we give a proof for completeness. The proof is essentially that of \cite{Klauder:2011}. Note that  our assumptions are chosen according to our needs rather than to obtain the most general statement. 
\begin{prop}[Berezin-Lieb inequalities]
Let $\Phi$ be a non-negative, concave continuous function on $[0,1]$, and let 
 $T$ be a positive trace class operator with $\tr(T)=1$. If $A$ is a positive compact operator on $L^2(\Rd)$ with $\|A\|_{\bo}\leq 1$, then
\begin{equation} \label{eq:berliebgen1}
    \int_{\Rdd} \Phi \circ (A\star T)(z) \ dz  \geq \tr(\Phi(A)),
\end{equation}
where $\Phi(A)$ is defined by the functional calculus. 
If $f\in L^1(\Rdd)$ is a positive function with $0\leq f \leq 1$, then 
\begin{equation} \label{eq:berliebgen2}
    \tr(\Phi(f\star T)) \geq \int_{\Rdd} \Phi \circ f (z) \ dz
\end{equation}
\end{prop}

\begin{proof}
 We can use the spectral theorem for self-adjoint compact operators to write $A$ in terms of its eigenvalues $\lambda_n$ and eigenfunctions $\varphi_n$ by \[A=\sum_{n=1}^\infty \lambda_n \varphi_n \otimes \varphi_n.\]
By the linearity of the convolutions, we find that 
\begin{align*}
    A\star T(z)&=\sum_{n=1}^\infty \lambda_n  (\varphi_n \otimes \varphi_n)\star T(z) \\
    &= \sum_{n=1}^\infty \lambda_n \langle T\pi(z)^* \varphi_n,\pi(z)^* \varphi_n  \rangle_{L^2}
\end{align*}
where the last equality follows from expanding the definition of the convolution $(\varphi_n \otimes \varphi_n) \star T$. Now note that for each fixed $z$, the sequence $\langle T\pi(z)^* \varphi_n,\pi(z)^* \varphi_n  \rangle_{L^2}$ is a probability distribution over $\mathbb{N}$: it is positive for each $n$ as $T$ is a positive operator, and the sequence sums to $1$ as $\{\pi(z)^*e_n\}_{n=1}^\infty$ is an orthonormal basis for each $z$ since $\pi(z)^*$ is unitary, so that
\[\sum_{n=1}^\infty \langle T\pi(z)^* \varphi_n,\pi(z)^* \varphi_n  \rangle_{L^2} = \tr(T)=1. \]

Hence we can apply Jensen's inequality for concave functions to get that 
\begin{align*}
    \Phi \circ (A\star T)(z)  &=   \Phi\left(\sum_{n=1}^\infty \lambda_n \langle T\pi(z)^* \varphi_n,\pi(z)^* \varphi_n  \rangle_{L^2} \right)  \\
    &\geq \sum_{n=1}^\infty  \Phi(\lambda_n) \langle T\pi(z)^* \varphi_n,\pi(z)^* \varphi_n  \rangle_{L^2}. 
\end{align*}
Note that we stay within the domain of $\Phi$: by Proposition \ref{prop:properties}, $|T\star A(z)|\leq \|A\|_{\bo} \|T\|_{\tco}=1$ for all $z$, and $\lambda_n \leq 1$ as $\|A\|_{\bo}\leq 1$.
Also note that by the same proposition 
\[\int_{\Rdd}   \langle T\pi(z)^* \varphi_n,\pi(z)^* \varphi_n  \rangle_{L^2}\ dz = \int_{\Rdd} (\varphi_n \otimes \varphi_n)\star T(z) \ dz = \tr(T)\|\varphi_n\|_{L^2}^2=1. \]
So by integrating the inequality above and changing the order of the sum and integral, we get
\begin{align*}
\int_{\Rdd} \Phi \circ (A\star T)(z) \ dz  &\geq \sum_{n=1}^\infty \Phi(\lambda_n) \int_{\Rdd}  \langle T\pi(z)^* \varphi_n,\pi(z)^* \varphi_n  \rangle_{L^2} \ dz \\
&= \sum_{n=1}^\infty \Phi(\lambda_n) \\
&= \tr(\Phi(A)).
\end{align*}
Turning to the second inequality, we will use the spectral theorem to write the positive trace class operator $f\star T$ in terms of its eigenvalues $\mu_n$ and eigenfunctions $\xi_n$ by 
\[f\star T=\sum_{n=1}^\infty \mu_n \xi_n \otimes \xi_n.\]
Since Proposition \ref{prop:properties} gives that $$\mu_n \leq \|f\star T\|_{\bo}\leq  \|f\|_{L^\infty} \|T\|_{\tco}\leq 1,$$ the eigenvalues $\mu_n$ belong to the domain of $\Phi$ and 
\[\Phi(f\star T)=\sum_{n=1}^\infty \Phi(\mu_n) \xi_n \otimes \xi_n,\]
and in particular 
\[\langle \Phi(f\star T)\xi_n,\xi_n \rangle_{L^2}=\Phi(\mu_n)=\Phi(\langle f\star T\xi_n,\xi_n \rangle_{L^2}).\]
This means that 
\begin{align*}
    \tr(\Phi(f\star T))&= \sum_{n=1}^\infty \langle \Phi(f\star T)\xi_n,\xi_n \rangle_{L^2} \\
    &= \sum_{n=1}^\infty \Phi(\langle f\star T\xi_n,\xi_n \rangle_{L^2}) \\
    &= \sum_{n=1}^\infty \Phi\left(\int_{\Rdd} f(z) \langle T\pi(z)^*\xi_n,\pi(z)^*\xi_n \rangle_{L^2} \ dz \right).
\end{align*}
Now note that for each fixed $n$, the function $\langle T\pi(z)^*\xi_n,\pi(z)^*\xi_n \rangle_{L^2}$ is a probability distribution over $\Rdd$: it is non-negative by the positivity of $T$ and $\int_{\Rdd} \langle T\pi(z)^*\xi_n,\pi(z)^*\xi_n \rangle_{L^2} \ dz=1$. So we can apply Jensen's inequality to get 
\begin{align*}
\tr(\Phi(f\star T))&=\sum_{n=1}^\infty \Phi\left(\int_{\Rdd} f(z) \langle T\pi(z)^*\xi_n,\pi(z)^*\xi_n \rangle_{L^2} \ dz \right)\\&\geq \sum_{n=1}^\infty \int_{\Rdd} \Phi\circ f(z) \langle T\pi(z)^*\xi_n,\pi(z)^*\xi_n \rangle_{L^2} \ dz  \\
&=  \int_{\Rdd} \Phi\circ f(z) \sum_{n=1}^\infty \langle T\pi(z)^*\xi_n,\pi(z)^*\xi_n \rangle_{L^2} \ dz \\
&= \int_{\Rdd} \Phi\circ f(z) \tr(T) \ dz \\
&= \int_{\Rdd} \Phi\circ f(z) \ dz.
\end{align*}
\end{proof}

When we pick $S$ to be a rank one operator $g\otimes g$, the associated mixed-state localization operator $\chi_\Omega\star (g\otimes g)$ is a kind of operator that has been studied extensively in the field of time-frequency analysis, namely a localization operator. By writing out the definition of the convolution, we find that 
\begin{equation*}
  (\chi_{\Omega}\star (g \otimes g))(f)=\int_{\Omega} V_{g}f(z) \pi(z) g \ dz.
\end{equation*}

\subsubsection{Entropy}
We have introduced the essential dimension of $S$ as the von Neumann entropy $H_{vN}(S)=-\tr(S\log (S))$. This definition applies the function $-x\log x$ to the operator $S$ using the so-called functional calculus, but an equivalent definition is that \begin{equation*}
    H_{vN}(S)=-\sum_{k} \lambda_k \log \lambda_k,
\end{equation*}
where $\lambda_k$ are the eigenvalues of $S$ from \eqref{Def:DatOpKL}. This is simply the Shannon entropy of the sequence of eigenvalues of $S$. 

As we will also be working with functions on $\Rdd$, we will also need the \textit{differential entropy} of a probability distribution $f$ on $\Rdd$, i.e. a positive function with $\int_{\Rdd} f(z) \, dz=1$. The differential entropy is then 
\begin{equation*}
    H(f)=-\int_{\Rdd} f(z) \log f(z) \, dz.
\end{equation*}

\subsubsection{Convolutional neural networks and data operator}
%...give motivation by stating a concrete version of the connection, references. 
Time-series such as audio data are routinely pre-processed before being used as input to CNNs. 
The  pre-processing steps  extract time-frequency structure 
of the signal 
which is  essential to human perception.  Thus, the pre-processing also
introduces invariance to
phase-shifts. Most commonly, the first processing step consists of taking a 
short-time Fourier transform, cf.~\eqref{Eq:STFT} followed by a non-linearity in the form of either
absolute value or absolute value squared. For an input signal $f$ one thus obtains
a first feature stage $F^0$ as given in \eqref{Eg:F0}. 
Due to the structure of the convolutional layers in a CNN, %as described in the introduction, 
the output of the first convolutional layer can then be written as 
 \begin{align*}
    F^1 (z,k) &= (F^0\ast m_k )(z) \\
    &=\langle \int_{z'} V_g f(z' )\cdot m_k (z-z' ) \pi (z' ) g\, dz' , f\rangle\\
    &=\langle(T_z \check{m}_k\star g\otimes g) f, f\rangle = \check{m}_k\ast [(f\otimes f)\star (\check{g}\otimes \check{g})](z)\\
    &=[\check{m}_k\star (f\otimes f)]\star (\check{g}\otimes \check{g})(z)
\end{align*} and these equalities offer various different interpretations of the first convolutional layer's output.  First, we can obviously interpret the value at
each $z$ as the correlation between the TF-localization operator 
$T_z \check{m}_k\star g\otimes g$ applied to $f$ with $f$. 
This operator itself, then, is simply the
TF-shifted version of $ \check{m}_k\star g\otimes g$.
More interestingly, we may equally consider
the operator $\check{m}_k\star (f\otimes f)$,
which now depends on the data points.
Ideally, the $m_k$ can be chosen as to
optimally enforce correlation within data
classes and to separate, or de-correlate, 
data points from different classes.

In order to isolate and investigate the question of correlations
between data points, we generalize the rank-one
operator $f\otimes f$ to the sum 
over all data points, 
thus considering $S = \sum_i f_i\otimes f_i$ and the
corresponding mixed-state localization operator 
$\chi_\Omega \star S$, which actually represents 
the local time-frequency averages of the 
data-points $f_i$. In CNNs, the lower layers 
realize convolutions with kernels $m_k$, where 
each kernel is supported in $\Omega$ and its
coefficients are learned from the data. The
intuitive task of the learnt kernels is to
strengthen correlation between data points from 
the same class, i.e. increase correlation within
the data points comprising $S$.

\subsubsection{Cohen's class}

We define the \textbf{Cohen's class associated with an operator} $S$ by 
\begin{equation*} 
  Q(f)=Q_S(f) = \check{S}\star (f\otimes f),\quad\text{for}~~f\in L^2(\mathbb{R}).
\end{equation*}

\begin{figure}
        \centering
        \includegraphics[width=1\textwidth]{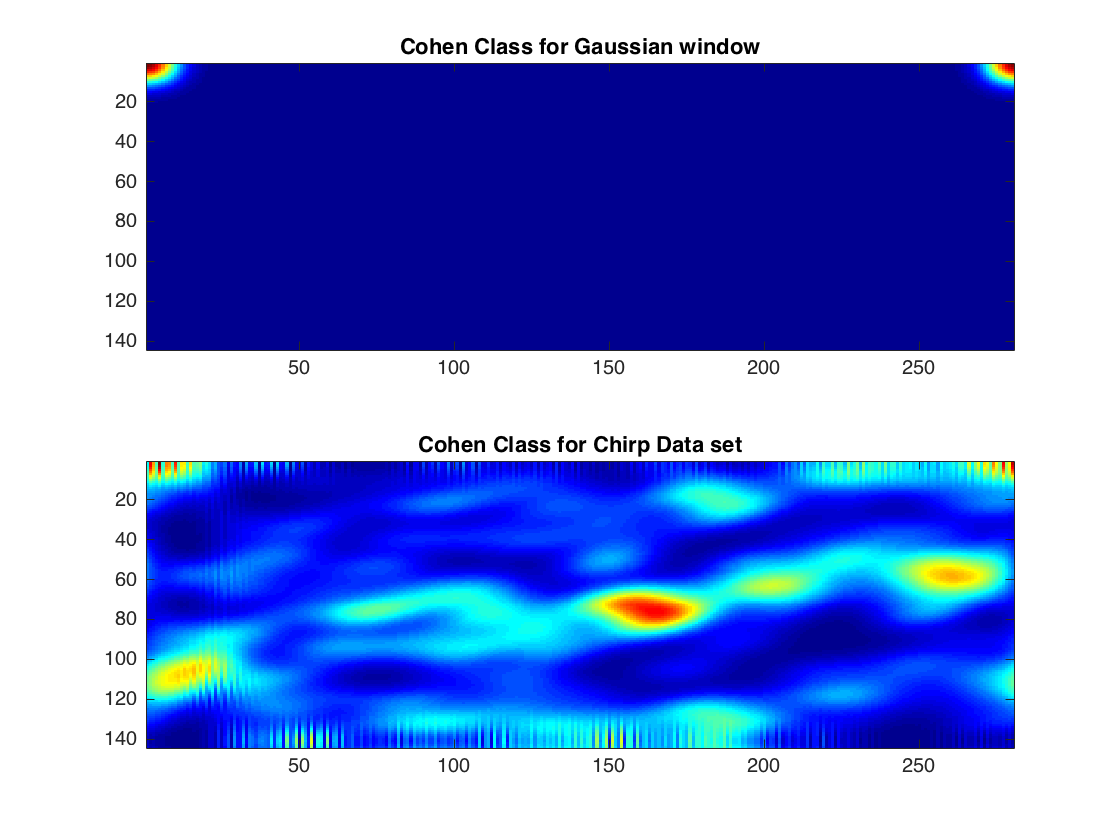}
        \caption{Cohen Class for $S = g\otimes g$, $g$ a Gaussian (upper plot) and for $S_{ch}$, the data operator of the chirp data set.}
        \label{fig: CohenCex}
    \end{figure}

An important example is the Cohen class associated with a rank one operator $g\otimes g$, which is the spectrogram
\begin{equation*}
  Q_{g \otimes g}(f)=(\check{g} \otimes \check{g}) \star (f\otimes f)=|V_{g}f|^2.
\end{equation*}
We will often meet operators $S$ that are linear combinations of rank one operators $S=\sum_k c_k g_k \otimes g_k$, which by a straightforward linearity property of the convolutions means that
\begin{equation*}
    Q_S(f)=\sum_k c_k |V_{g_k}f|^2.
\end{equation*}
\begin{rem}
Cohen's class distributions were introduced by Cohen in the context of quantum 
mechanics, cf.~\cite{Cohen89}, and are defined by applying smoothing 
kernels to the so-called Wigner distribution.  For details on the connection between the two distinct definitions, see~\cite{LuSk20}. %nterpretation, simple examples,...
We give two examples of Cohen class analysis in Figure~\ref{fig: CohenCex}. In the upper plot, we show $Q_S(g)$ for 
$S = g\otimes g$, where $g$ a Gaussian window. The lower plot shows 
$Q_{S_{ch}}(g)$. 
\end{rem}

\subsubsection{Mixed-state localization operators, Cohen's class and data operators} 
We will now describe some connections between the  concepts introduced above. We begin by connecting the Cohen class $Q_S$ to the mixed-state localization operator (or $\Omega$-augmentation) $\frac{\chi_\Omega}{|\Omega|}\star S$.

Since all convolutions preserve positivity, i.e.\,the convolution of a positive function with a positive operator is a positive operator, we get whenever $S$ is a positive trace class operator with $\tr(S)=1$ and $\Omega$ is a compact domain that $\chi_\Omega \star S$ is a positive, compact operator. Hence  $\chi_\Omega \star S$ has a diagonalization
\begin{equation*}
    \chi_\Omega \star S =\sum_{k=1}^\infty \lambda_k^\Omega h_k^\Omega \otimes h_k^\Omega, 
\end{equation*}
where $0\leq \lambda_k^\Omega\leq 1$ are the eigenvalues of $\chi_\Omega \star S$ with eigenfunctions $h_k^\Omega$ and $\sum_{k=1}^\infty \lambda_k^\Omega = \tr(\chi_\Omega \star S)=|\Omega|$.

In particular, $\chi_\Omega \star S$ is completely determined by its eigenvectors and eigenvalues, which are further determined by $\Omega$ and the Cohen class $Q_S$:

  \begin{align*}
  \lambda_k^\Omega &= \max_{\|f\|_2=1, f \perp h_1^\Omega, \dots h_{k-1}^\Omega} \int_{\Omega} Q_{S}(f)(z) \ dz \\
      h_k^\Omega &= \argmax_{\|f\|_2=1, f \perp h_1^\Omega, \dots h_{k-1}^\Omega} \int_{\Omega} Q_{S}(f)(z) \ dz.
  \end{align*}

By Lemma 3.1 in \cite{Werner} we have whenever $\tr(S)=1$ and $\|f\|_2=1$ that 
\begin{equation*}
    \int_{\Rdd} Q_S(f)(z) \ dz =1.
\end{equation*}

Since $S$ is also a positive operator, it is easy to show that $Q_S(f)$ is a positive function. We can therefore measure the localization of $f$ in $\Omega$ by how much of $Q_S(f)$ is concentrated in $\Omega$, i.e. by $\int_{\Omega} Q_{S}(f)(z) \ dz$. If we use this notion of localization, the expressions above say that the eigenfunctions $h_k^\Omega$ form a maximally localized orthonormal system. Of course, this notion of localization depends heavily on $S$, and may differ profoundly from more intuitive notions of time-frequency localization.

The expression for $\lambda_k^\Omega$ shows that $\lambda_k^\Omega$ is a measure of the localization of the maximally localized orthonormal system $h_k^\Omega$. The size of $\lambda_k^\Omega$ for varying $k$ will therefore depend on the existence of well-localized functions $h_k^\Omega$, as measured by integrating $Q_S(h_k^\Omega)$ over $\Omega$. Since the entropy of $\frac{\chi_\Omega}{|\Omega|}\star S$ can be expressed in terms of the Shannon entropy of $\frac{\lambda_k^\Omega}{|\Omega|}$, the entropy of $\frac{\chi_\Omega}{|\Omega|}\star S$ is also related to the existence of such well-localized functions. 

\subsubsection{Data Operator}\label{Sec:DatOp} 
So far we have worked with positive trace class operators $S$, sometimes normalized by $\tr(S)=1$. As we have already hinted at, we will restrict our attention to those $S$ arising as data operators as in \eqref{Def:DatOp}. This means that 
we let $S=\sum_{i=1}^N f_i\otimes f_i$ with normalized data $\sum_{i=1}^N \|f_i\|_2^2=1$. The associated Cohen's class distribution is then 
\begin{equation*}
   Q_S(f)(z)=\sum_{i=1}^N |V_{f_i}f(z)|^2, 
\end{equation*}
 and the total correlation function $\tilde{S}=S\star \check{S}$ is given by
\begin{equation*}
    \tilde{S}(z)=\sum_{i,j=1}^N |V_{f_i}f_j(z)|^2= \sum_{i,j=1}^N |\langle f_j,\pi(z) f_i \rangle|^2.
\end{equation*}
This expression is the background for the terminology of total correlation function: it measures the total sum of the correlations (as measured by inner products) between time-frequency shifted data points.
\begin{rem}
For $\{f_i=\pi(\lambda_i)g:\,i=1,...,N\}$ the total correlation function $\tilde{S}(z)=\sum_{i,j=1}^N |\langle g,\pi(z+\lambda_i-\lambda_j) g \rangle|^2$. Note that $\tilde{S}(0)=\sum_{i,j=1}^N|\langle g,\pi(\lambda_i-\lambda_j) g \rangle|^2$ is Hilbert-Schmidt norm of the Gramian of the system $\{\pi(\lambda_i)g:\,i=1,...,N\}$.
\end{rem}
For this choice of $S$ our $\mathrm{argmax}$-expression for the eigenfunctions of $\frac{\chi_\Omega}{|\Omega|}\star S$ also gets a simplified expression:
\begin{equation*}
      h_k^\Omega = \argmax_{\|f\|_2=1, f \perp h_1^\Omega, \dots h_{k-1}^\Omega} \sum_{i=1}^N \int_{\Omega} |V_{f_i}f(z)|^2 \ dz,
  \end{equation*}
which can be interpreted as saying that $h_k^\Omega$ are (locally in $\Omega$) maximally correlated with the data points $f_i$.

\subsection{Correlations and essential dimensions} 
From now on, $S$ will denote a data operator $S=\sum_{i=1}^N f_i\otimes f_i$ as above.
Our main goal is to understand the essential dimension of the data operator $S$ and its $\Omega$-augmentation $\frac{\chi_\Omega}{|\Omega|}\star S$, and how this relates to the total correlation function $\tilde{S}$. We propose that all of these give different perspectives on correlations in our data set. Since essential dimension is defined in terms of von Neumann entropy, we first investigate how entropy fits our framework. 

\subsubsection{Entropy in quantum harmonic analysis}
We recall that the von Neumann entropy $H_{vN}(S)$ of a positive trace class operator $S$  with $\tr(S)=1$ is given by 
\begin{equation*}
    H_{vN}(S)=-\tr(S\log S)=-\sum_{k=1}^\infty \lambda_k \log \lambda_k
\end{equation*}
where $\lambda_k$ denotes the eigenvalues of $S$. In the 
case of a data operator $S_\D$, this entropy reflects the
effective dimensionality of the data set $\D$.  Since the function $t\mapsto -t\log t$ is concave on $[0,1]$, we can use the Berezin-Lieb inequalities due to Werner \cite{Werner} to obtain the following  relations. 

\begin{prop}[Berezin-Lieb inequalities for entropy] \label{prop:berlieb}

Let $\Omega \subset \Rdd$ be a compact subset of $\Rdd$. Then, with $H$ denoting the differential entropy, 
 \begin{equation} \label{eq:berlieb1}
     H\left(\frac{\chi_\Omega}{|\Omega|}\ast \tilde{S}\right)   \geq H_{vN}\left( \frac{\chi_\Omega}{|\Omega|}\star S \right)\geq \log |\Omega|.
    \end{equation}
    \begin{equation}\label{eq:berlieb2}
        H(\tilde{S})\geq H_{vN}(S).
    \end{equation}
\end{prop}

\begin{proof}
Throughout the proof we will refer to \eqref{eq:berliebgen1} and $\eqref{eq:berliebgen2}$ with $\Phi(x)=-x\log x$.

We start with the relation 
\[H_{vN}\left( \frac{\chi_\Omega}{|\Omega|}\star S \right)\geq \log |\Omega|.\]
This is trivially true if $|\Omega|<1$, as the left hand side is non-negative and the right hand side is negative in that case. So we may assume that $|\Omega|\geq 1$, and the relation then follows by \eqref{eq:berliebgen2} with $f=\frac{\chi_\Omega}{|\Omega|}$ and $T=S$.

The other part of  \eqref{eq:berlieb1} follows by picking $T=\check{S}$ and $A=\frac{\chi_\Omega}{|\Omega|}\star S$ in \eqref{eq:berliebgen1}. This leads to \[H\left( \left(\frac{\chi_\Omega}{|\Omega|}\star S\right) \star  \check{S}\right)   \geq H_{vN}\left( \frac{\chi_\Omega}{|\Omega|}\star S \right),\] and we can use associativity of convolutions to get 
\[H\left( \left(\frac{\chi_\Omega}{|\Omega|}\star S\right) \star  \check{S}\right) =H\left( \frac{\chi_\Omega}{|\Omega|}\ast (S \star  \check{S})\right)=H\left( \frac{\chi_\Omega}{|\Omega|}\ast \tilde{S}\right). \]
Finally, to get \eqref{eq:berlieb2} simply use $A=S$ and $T=\check{S}$ in \eqref{eq:berliebgen1}.
\end{proof}  

These inequalities express the fact that convolutions increase the entropy, even when we use several different definitions of convolutions and both von Neumann and differential entropy. In particular, we get upper bounds for the essential entropy of both $S$ and its $\Omega$-augmentation in terms of the \textit{differential entropy} of the total correlation function. 

It is easy to check that the total correlation function $\tilde{S}$ is non-negative with $\int_{\Rdd} \tilde{S}(z) \ dz=1$, $ \tilde{S}(z)\leq \tilde{S}(0)=\tr(S^2)\leq 1$
and $\tr(S^2)=1$ if and only if $S$ is a rank one operator. We may interpret the entropy $H(\tilde{S})$ as 
a measure of TF-concentration of $\tilde{S}$. 
In other words: how large is the subset of $\R^2$ such that $\tilde{S}(z)$ is non-negligible? A small value for $H(\tilde{S})$ suggests that $\tilde{S}$ is only non-negligible in a small region around $z=0.$ Since 
\begin{equation*}
    \tilde{S}(z)=\sum_{i,j=1}^N |V_{f_i}f_j(z)|^2= \sum_{i,j=1}^N |\langle f_j,\pi(z) f_i \rangle|^2,
\end{equation*}
the correlations in the data are then captured by \textit{small} time-frequency shifts $\pi(z)f_i$ of the data points. % \todo{this is vague... We need to decide when and how to bring in the correlations term}. 
By equation \eqref{eq:berlieb2} in  Proposition \ref{prop:berlieb}  a thus well-localized total correlation function implies a small essential dimension of the data set. This observation underlines the connection 
between local correlations and effective dimensionality.

The connection between the concentration of $\tilde{S}$ and the essential dimension is also captured if we measure the concentration of $\tilde{S}$ by $\int_{\Rdd} |z|^2 \tilde{S}(z) \ dz$, as the following proposition shows.

\begin{prop}
    Let $f\in L^1(\R^d)$ be non-negative with $\int_{\Rd} f(t) \ dt=1.$ Then
    \begin{equation*}
        e^{H(f)/d} \leq \sqrt{\frac{2\pi e}{d}} \sqrt{\int_{\Rd} |x|^2 f(x) \ dx-|\mu|^2,}
    \end{equation*}
    where $\mu \in \mathbb{R}^d$ is the mean, i.e. $\mu_i=\int_{\Rd} x_i f(x) \ dx.$ In particular if $S$ is a positive trace class operator with $\tr(S)=1$, then
     \begin{equation*}
        e^{H(\tilde{S})/2d} \leq \sqrt{\frac{\pi e}{d}} \sqrt{\int_{\Rdd} |z|^2 \tilde{S}(z) \ dz-|\mu_S|^2}
    \end{equation*}
    where $\mu_S$ is the mean of $\tilde{S}$.
\end{prop}
\begin{proof}
    By assumption, $f$ is a probability density function. We let $\Sigma$ be the associated covariance matrix. This matrix satisfies that
    \begin{equation} \label{eq:amgm}
        \frac{\tr(\Sigma)}{d}\geq \det (\Sigma)^{1/d}.
    \end{equation} 
    To see this, recall that $\Sigma$ is a positive semi-definite matrix. If $\lambda_i$ are its (non-negative) eigenvalues, we know that $\tr(\Sigma)=\sum_{i=1}^{d} \lambda_i$ and $\det(\Sigma)=\prod_{i=1}^d \lambda_i$. The inequality therefore states that
    $$\frac{\sum_{i=1}^{d} \lambda_i}{d}\geq \left(\prod_{i=1}^d \lambda_i\right)^{1/d},$$
    which is true by the AM-GM inequality. 
    
    We then recall an inequality due to Shannon, see \cite{fosi97} for a proof, namely that 
    \begin{equation} \label{eq:shannon}
        H(f)\leq \frac{1}{2} \log [(2\pi e)^d \det (\Sigma)],
    \end{equation}
    where $H(f)$ is the differential entropy.
    
    Let us then combine these two inequalities. By taking logarithms of \eqref{eq:amgm} we find that $$\log \det (\Sigma)\leq d \log \frac{\tr (\Sigma)}{d}.$$
    By inserting this into \eqref{eq:shannon} we get that
    \begin{align*}
        H(f) &\leq \frac{d}{2} \log (2\pi e)+\frac{1}{2} \log (\det (\Sigma)) \\
        &\leq  \frac{d}{2} \log (2\pi e)+\frac{1}{2} d \log \frac{\tr (\Sigma)}{d} \\
        &= \frac{d}{2} \log [2\pi e \tr(\Sigma)/d].
    \end{align*}
    Hence we end up with
    $$e^{H(f)/d} \leq \sqrt{2\pi e \tr(\Sigma)/d}=\sqrt{2\pi e/d} \sqrt{\tr(\Sigma)}.$$
    The desired inequality follows, as $$\tr(\Sigma)=\int_{\Rd} |x|^2 f(x) \ dx-|\mu|^2$$ by the definition of the covariance matrix. In particular, the inequality for $\tilde{S}$ holds as $f=\tilde{S}$ satisfies the assumptions of the first part. 
\end{proof}

\begin{rem}
We also mention that Huber et al. have given an improvement on the relation between the essential dimensions of $S$ and $f\star S$ where $f$ is some probability distribution on $\Rdd$, namely the following \textit{entropy power inequality}\cite{Huber:2017}:
\begin{equation*}
    e^{H_{vN}(f \star S)/d}\geq e^{H(f)/d}+e^{H_{vN}(S)/d}.
\end{equation*}

The precise assumptions on $f$ and $S$ for this to hold are unfortunately not clear in \cite{Huber:2017}. 
\end{rem}

\subsubsection{The projection functional and average lack of concentration}
From the Berezin-Lieb inequalities, we know that the spread of the total correlation function, as measured by its differential entropy, gives an upper bound for the essential dimension. When we add $\Omega$ to the setup and consider the $\Omega$-augmentation $\frac{\chi_\Omega}{|\Omega|}\star S$, we also have the upper bound 
\begin{equation*}
     H\left(\frac{\chi_\Omega}{|\Omega|}\ast \tilde{S}\right)   \geq H_{vN}\left( \frac{\chi_\Omega}{|\Omega|}\star S \right),
\end{equation*}
which gives a first indication that the essential dimension of the $\Omega$-augmentation depends on the interplay of $\Omega$ and $\tilde{S}$. Our goal is now to understand this interplay in detail.

The first quantity we will use to study this, is the \textit{projection functional}
\begin{equation*}
    P(\chi_\Omega \star S)=\tr \left( \chi_\Omega \star S \right) - \tr \left( (\chi_\Omega \star S)^2 \right).
\end{equation*}
Clearly, the projection functional measures how much $\chi_\Omega \star S$ deviates from being a projection. When $\chi_\Omega \star S$ is a projection, we are in the idealized situation where $|\Omega|$ is an integer and the eigenvalues are $\lambda_k^\Omega=1$ for $k\leq |\Omega|$ and $\lambda_k^\Omega=0$ otherwise. This would mean that the $\Omega$-augmentation is completely described by $|\Omega|$ eigenfunctions:
\begin{equation} \label{eq:locopisprojection}
    \chi_\Omega \star S=\sum_{k=1}^{|\Omega|} h_k^\Omega \otimes h_k^\Omega.
\end{equation}
In addition, it is easy to show that this situation gives $H_{vN}\left(\frac{\chi_\Omega}{|\Omega|}\star S\right)=\log |\Omega|$, which by Proposition \ref{prop:berlieb} means that we are in the situation of minimal essential dimension of the $\Omega$-augmentation.

In fact, since 
\begin{equation} \label{eq:projvsent}
x-x^2\leq -x\log x \text{ for } x\in [0,1]
\end{equation} the projection functional is clearly related to the essential dimension.

\begin{lem} \label{lem:lowerbound}
For $\Omega\subset \Rdd$ compact and $S$ a positive trace class operator with $\tr(S)=1$, 
\begin{equation*}
    \log |\Omega|+\frac{1}{|\Omega|} P(\chi_\Omega \star S)\leq H_{vN}\left( \frac{\chi_\Omega}{|\Omega|}\star S \right).
\end{equation*}
\end{lem}
\begin{proof}
It is straightforward to show that $P(\chi_\Omega \star S)=\sum_{k=1}^\infty\left( \lambda_k^\Omega - \left(\lambda_k^\Omega\right)^2\right)$. Since $0\leq \lambda_k^\Omega \leq 1$, \eqref{eq:projvsent} gives that 
\begin{equation*}
    P(\chi_\Omega \star S)\leq \sum_{k=1}^\infty -\lambda_k^\Omega \log \lambda_k^\Omega.
\end{equation*}
Then note that 
\begin{align*}
    H\left( \frac{\chi_\Omega}{|\Omega|}\star S \right)&=\sum_{k=1}^\infty -\frac{\lambda_k^\Omega}{|\Omega|} \log \frac{\lambda_k^\Omega}{|\Omega|} \\
    &= \log |\Omega| \sum_{k=1}^\infty  \frac{\lambda_k^\Omega}{|\Omega|}-\frac{1}{|\Omega|}\sum_{k=1}^\infty \lambda_k^\Omega \log \lambda_k^\Omega \\
    &= \log |\Omega|-\frac{1}{|\Omega|}\sum_{k=1}^\infty \lambda_k^\Omega \log \lambda_k^\Omega \\
    &\geq \log |\Omega|+\frac{1}{|\Omega|}\cdot P(\chi_\Omega \star S). \qedhere
\end{align*}
\end{proof}

We are interested in knowing what properties of $\Omega$ and the data set make $P(\chi_\Omega \star S)$ small or large. In light of the above, this means that we want to know what prevents the essential dimension of $\frac{\chi_\Omega}{|\Omega|}\star S$ from being small, or by \eqref{eq:locopisprojection} what prevents the $\Omega$-augmentation to be completely described by a few functions. To clarify this, we use the following result which is an easy consequence of \cite[Lem. 4.2]{LuSk20}.

\begin{prop} \label{prop:projfunisalc}
For $\Omega\subset \Rdd$ compact and $S$ a positive trace class operator with $\tr(S)=1$, 
\begin{equation*}
    P(\chi_\Omega \star S)=\int_{\Omega} \left( 1- \int_{\Omega-z} \tilde{S}(z') \, dz' \right) \, dz.
\end{equation*}
\end{prop}
To help understand this result, we introduce the \textit{average lack of concentration of $\tilde{S}$ with respect to $\Omega$} as
\begin{equation} \label{eq:definealc}
    \alc(\tilde{S},\Omega) := \frac{1}{|\Omega|} \int_{\Omega} \left( 1- \int_{\Omega-z} \tilde{S}(z') \, dz' \right) \, dz.
\end{equation} 
The reader will of course have noticed that $\alc(\tilde{S},\Omega)=\frac{1}{|\Omega|}P(\chi_\Omega \star S)$ by Proposition \ref{prop:projfunisalc}. Nevertheless, we single out $\alc(\tilde{S},\Omega)$ in order to discuss its interpretation. We start by considering the integrand in \eqref{eq:definealc} for fixed $z$:
\begin{equation} \label{eq:lc}
    1- \int_{\Omega-z} \tilde{S}(z') \, dz'.
\end{equation}
As we know that $\int_{\Rdd} \tilde{S}(z') \ dz'=1$, \eqref{eq:lc} measures how far $\tilde{S}$ is from being completely concentrated in $\Omega-z$ --- i.e. the lack of concentration. The average lack of concentration $\alc(\tilde{S},\Omega)$ is obtained by averaging the lack of concentration in \eqref{eq:lc} over $\Omega$. In this sense, $\alc(\tilde{S},\Omega)$ measures how far $\tilde{S}$ is from being concentrated in $\Omega-z$ as $z$ varies over $\Omega$. We can restate Lemma \ref{lem:lowerbound} and Proposition \ref{prop:berlieb} to get upper and lower bounds for the essential dimension of the $\Omega$-augmentation in terms of the function $\tilde{S}$ and the domain $\Omega$. 

\begin{thm}\label{thm:ALC_ED}  
For $\Omega\subset \Rdd$ compact and $S$ a positive trace class operator with $\tr(S)=1$, 
\begin{equation*}
    \log |\Omega|+ \alc(\tilde{S},\Omega) \leq  H_{vN}\left( \frac{\chi_\Omega}{|\Omega|}\star S \right) \leq H\left(\frac{\chi_\Omega}{|\Omega|} \ast \tilde{S} \right).
\end{equation*}
\end{thm}
If the data set consists just of one point $f\otimes f$, then $\frac{\chi_\Omega}{|\Omega|}\star S$ is the localization operator $A_\Omega^f\xi=\int_\Omega V_f\xi(z)\pi(z)\xi\,dz$, $\tilde{S}=|V_ff|^2$ and $\alc(\tilde{S},\Omega)$ is the projection functional. Thus the preceding result is a statement about the von Neumann entropy of the localization operator $A_\Omega^f$ and the spectrogram of the data point $|V_ff|^2$.

We see that the essential dimensionality of the $\Omega$-augmentation really depends on $\alc(\tilde{S},\Omega)$, i.e. on how concentrated $\tilde{S}$ is in the sets $\Omega-z$ for $z\in \Omega$. We emphasize that $\alc(\tilde{S},\Omega)$ highly depends on the interplay between $\Omega$ and $\tilde{S}$, and that this interplay seems to be the key to understanding the essential dimension of the $\Omega$-augmentation. However, it is possible to get an upper bound that decouples the effects of $\Omega$ and $\tilde{S}$, see the proof of \cite[Lem. 5.3]{LuSk20}.

\begin{prop}
For $\Omega\subset \Rdd$ compact and $S$ a positive trace class operator with $\tr(S)=1$, 
\begin{equation*}
    \alc(\tilde{S},\Omega)\leq  \frac{ |\partial \Omega|}{|\Omega|}  \int_{\Rdd} \tilde{S}(z) |z| \ dz.
\end{equation*}
\end{prop}
\begin{rem}
In this result, the size $|\partial \Omega|$ of the perimeter is defined as the variation of $\chi_\Omega$, in other words
\begin{equation*}
  |\partial \Omega|=\sup \left\{ \int_{\Rd} \chi_\Omega(x) \text{div} \phi (x) \ dx : \phi \in C_c^1(\Rd,\Rd), |\phi(x)|\leq 1\ \forall x\in \Rd  \right\},
\end{equation*}
where $\text{div}\phi$ is the divergence of $\phi$,
$C_c^1(\Rd,\Rd)$ is the set of compactly supported
differentiable functions from $\Rd$ to $\Rd$ and $|\phi(x)|$ denotes the Euclidean norm on $\Rd$. This is the same measure of the perimeter used in \cite{LuSk20,Abreu:2016,Abreu:2017}.
\end{rem}

\subsection{More relations on ALC
and approximation of data operators}

We will also mention how $\alc(\tilde{S},\Omega)$ influences 
other quantities we have discussed. We begin by bounding how 
much of the $\Omega$-augmentation is captured by the first few eigenfunctions. 
For this, we first need a simple lemma; its proof can be 
deduced from the proof of \cite[Thm. 6.1]{LuSk20}.

\begin{lem}
For $\Omega\subset \Rdd$ compact and $S$ a positive trace class operator with $\tr(S)=1$, and let $A_\Omega=\lceil |\Omega| \rceil$. Then
\[
    \alc(\tilde{S},\Omega)\geq 1-\sum_{k=1}^{A_\Omega} \frac{\lambda_k^\Omega}{|\Omega|}.
\]
\end{lem}

\begin{proof}
By Proposition \ref{prop:projfunisalc} we have that 
\begin{align*}
    |\Omega|\alc(\tilde{S},\Omega)&=P(\chi_\Omega \star S) \\
    &= \tr(\chi_\Omega \star S)-\tr((\chi_\Omega \star S)^2) \\
    &= \sum_{k=1}^{A_\Omega} \lambda_k^\Omega (1-\lambda_k^\Omega)+\sum_{k=A_\Omega+1}^\infty  \lambda_k^\Omega (1-\lambda_k^\Omega) \\
		&\geq \lambda_{A_\Omega}^\Omega  \sum_{k=1}^{A_\Omega} (1-\lambda_k^\Omega) +(1-\lambda_{A_\Omega}^\Omega)\sum_{k=A_\Omega+1}^\infty  \lambda_k^\Omega \\
		&= \lambda_{A_\Omega}^\Omega A_\Omega - \lambda_{A_\Omega}^\Omega |\Omega|+\sum_{k=A_\Omega +1}^{\infty} \lambda_k^\Omega \\
		&= \lambda_{A_\Omega}^\Omega(A_\Omega-|\Omega|) + |\Omega|-\sum_{k=1}^{A_\Omega}\lambda_k^\Omega \\ 
		&\geq |\Omega|-\sum_{k=1}^{A_\Omega}\lambda_k^\Omega. \qedhere
\end{align*}
\end{proof}
\begin{thm} \label{thm:finiterankapprox}
For $\Omega\subset \Rdd$ compact and $S$ a positive trace class operator with $\tr(S)=1$, let $T_\Omega=\sum_{k=1}^{A_\Omega} h_k^\Omega\otimes h_k^\Omega$ where $A_\Omega=\lceil |\Omega| \rceil$. Then
\begin{equation*}
    \frac{\|\chi_\Omega\star S-T_\Omega\|_{\tco}}{|\Omega|} \leq \frac{A_\Omega-|\Omega|}{|\Omega|}+ 2\cdot \alc(\tilde{S},\Omega).
\end{equation*}
\end{thm}
\begin{proof}
Using the singular value decomposition 
$\chi_\Omega \star S=\sum_{k=1}^\infty \lambda_k h_k^\Omega \otimes h_k^\Omega$ 
we find that 
\begin{align*}
	\left\|\chi_\Omega \star S- \sfr \right\|_{\tco}&= \left\|\sum_{k=1}^{\infty}\lambda_k^\Omega h_k^\Omega \otimes h_k^\Omega- \sum_{k=1}^{A_{\Omega}}h_k^\Omega \otimes h_k^\Omega\right\|_{\tco} \\
	&= \left\|\sum_{k=A_\Omega+1}^{\infty} \lambda_k^\Omega h_k^\Omega \otimes h_k^\Omega- \sum_{k=1}^{A_{\Omega}}(1-\lambda_k^\Omega)h_k^\Omega \otimes h_k^\Omega\right\|_{\tco} \\
	&= \sum_{k=1}^{A_\Omega} (1-\lambda_k^\Omega) + \sum_{k=A_\Omega +1}^\infty \lambda_k^\Omega.
\end{align*}
We use $|\Omega|=\sum_{k=1}^\infty \lambda_k^\Omega$ to find that 
\begin{align*}
    \sum_{k=1}^{A_\Omega} (1-\lambda_k^\Omega) + \sum_{k=A_\Omega +1}^\infty \lambda_k^\Omega&= A_\Omega -\sum_{k=1}^{A_\Omega} \lambda_k^\Omega+ |\Omega| - \sum_{k=1}^{A_\Omega} \lambda_k^\Omega \\
    &=  A_\Omega + |\Omega| - 2\sum_{k=1}^{A_\Omega} \lambda_k^\Omega \\
    &= A_\Omega-|\Omega| + 2\left( |\Omega| - \sum_{k=1}^{A_\Omega} \lambda_k^\Omega \right).
\end{align*}
The claim follows from the lemma.
\end{proof}

Asymptotically, as the size of $\Omega$ increases, it is clear that the average lack of concentration with respect to $\Omega$ will decrease. This is made formal in the next statement.

\begin{prop}
For $\Omega\subset \Rdd$ compact and $S$ a positive trace class operator with $\tr(S)=1$, with respect to  $R \Omega$, we have
\begin{equation*}
    \lim_{R\to \infty} \alc(\tilde{S},R\Omega)=0.
\end{equation*}
\end{prop}

\begin{proof}
An easy calculation gives that
\begin{equation*}
    \alc(\tilde{S},R\Omega)=1-\frac{1}{|R\Omega|} \int_{R\Omega} \int_{R\Omega} \tilde{S}(z-z')\ dz dz',
\end{equation*}
and the result is immediate from \cite[Cor. 3.4.1]{LuSk20}.
\end{proof}

\bibliographystyle{abbrv}{}
\bibliography{DimBib}   
\end{document}